\documentclass[11pt]{article}      
\usepackage{amsmath,amsthm,amsfonts,graphicx}
\usepackage{multirow}
\usepackage{multirow}
\usepackage{xcolor}
\def\smallddots{\mathinner{\raise7pt\hbox{.}\raise4pt\hbox{.}\raise1pt\hbox{.}}} 
\def\smallsdots{\mathinner{\raise1pt\hbox{.}\raise4pt\hbox{.}\raise7pt\hbox{.}}}

\DeclareMathOperator{\diag}{diag}
\DeclareMathOperator{\prob}{Probability}
\DeclareMathOperator{\rank}{rank}

\DeclareMathOperator{\cond}{cond}

\newtheorem{theorem}{Theorem}[section]
\newtheorem{outline}{Outline}[section]
\numberwithin{equation}{section}
\numberwithin{table}{section}
\newtheorem{lemma}{Lemma}[section]

\newtheorem{corollary}{Corollary}[section]

\newtheorem{algorithm}{Algorithm}[section]
\newtheorem{example}{Example}[section]
\newtheorem{definition}{Definition}[section]
\newtheorem{assumption}{Assumption}[section]

\newtheorem{remark}{Remark}[section]

   \usepackage{enumitem}
\newlist{steps}{enumerate}{1}
\setlist[steps, 1]{label = Step \arabic*:}
\usepackage{booktabs}
\usepackage{subcaption}
\usepackage{tikz}
\usetikzlibrary{patterns, patterns.meta}
\usepgflibrary[patterns]

\begin{document} 
 
\title{Novel Insights into  
Low Rank
Approximation 
at Sublinear Cost} 
\author{Victor Y. Pan}
\author{Soo Go$^{[2],[a]}$,
Qi Luan$^{[2],[b]}$, Victor Y. Pan$^{[1, 2],[c]}$, 
John Svadlenka$^{[2],[d]}$, Liang Zhao$^{[1, 2],[e]}$
\\\\
$^{[1]}$ Department of Computer Science \\
Lehman College of the City University of New York \\
Bronx, NY 10468 USA \\
$^{[2]}$ Ph.D. Programs in  Computer Science and Mathematics \\
The Graduate Center of the City University of New York \\
New York, NY 10016 USA \\
$^{[a]}$ sgo@gradcenter.cuny.edu\\
$^{[b]}$ qi\_luan@yahoo.com \\ 
$^{[c]}$ victor.pan@lehman.cuny.edu \\ 
http://comet.lehman.cuny.edu/vpan/  \\
$^{[d]}$jsvadlenka@gradcenter.cuny.edu \\ 
$^{[e]}$ Liang.Zhao1@lehman.cuny.edu \\
} 
\date{}

\maketitle

   
\begin{abstract}    
Low rank approximation of a matrix ({\em LRA}) is a highly important area of  Numerical Linear and Multilinear Algebra and  Data Mining and Analysis. One can operate  with an LRA superfast -- by using much fewer memory cells and flops than an input matrix  has entries.\footnote{``Flop"  stands for ``floating point  operation".} Can we, however, compute an LRA of a matrix superfast? YES and NO. For worst case inputs,
any LRA algorithm
fails miserably unless it involves all input entries, but in computational
practice worst case inputs seem to appear rarely, and accurate LRA are routinely computed  superfast for large and important classes of matrices,
in particular in the
memory efficient  form of CUR, widely used in data analysis. 
We advance formal study of this YES and NO coexistence
by proving  {\em novel universal} 
upper bounds   on the spectral output error norms of {\em all CUR LRA} algorithms and,
  under a fixed probabilistic structure in the space of input matrices, on
 both spectral and Frobenius  error norms  of {\em nearly  all sketching LRA algorithms}. These bounds imply
that superfast LRA algorithms
of the two  kinds fail miserably   only for a very narrow  input 
 class. Furthermore, in our numerical tests such  superfast
 algorithms 
were consistently
 much more accurate than our upper estimates ensure and usually were reasonably close to  
optimal.
\end{abstract}
 




\bigskip
\medskip


  

\section{Introduction}\label{sintr} 


{\bf 1.1.  Superfast LRA 
CUR LRA:  background.}
LRA of a matrix is among the most fundamental problems of 
Numerical Linear and Multilinear Algebra and Data  Mining and Analysis, with  
  applications ranging from machine
learning theory and neural networks
to term document data and DNA SNP data (see surveys \cite{HMT11,MT20,TWa}). 
For example, matrices that represent 
  Big Data (e.g., unfolding matrices of multidimensional tensors)
tend to be so immense that  only a tiny fraction of their entries
fits primary memory of a computer, but  
   quite typically  they admit LRA,
  that is,
can be closely approximated
with the products of a few matrices of smaller sizes \cite{UT19}. Such products
can be  multiplied by a vector
 {\em superfast}, aka at {\em sublinear cost}, -- by using much fewer memory cells and flops than the matrix has 
entries.\footnote{Superfast algorithms have been studied long and  extensively for
 structured matrices having small displacement rank  
(see more in Sec. 1.4.1), and we extend this concept to computations of LRA.} This implies dramatic acceleration of various   popular matrix computations  reduced to that basic operation, but
can we compute an LRA of a matrix superfast or at least {\em quasi superfast}, that is, by involving only a small fraction of all entries of an input matrix? YES and NO.

 NO, because on 
the worst case inputs and even on the small matrix families of our Example \ref{exdlt}  
 any  algorithm  {\em fails
miserably}
for computing accurate LRA and even for a posteriori  estimation
of its accuracy
 unless it involves
all entries of
an  input matrix.

YES, because in computational
practice worst case inputs seem to appear rarely, and 
 accurate  LRAs are
 routinely  computed  superfast for a large and important  class of
matrices,
in  particular  by means of the popular 
{\em Cross-Approximation (C--A)} iterations, specializing  the Alternating Directions Implicit (ADI) method \cite{S16} to LRA and recalled in Sec. \ref{scait}.

Their output LRAs are memory efficient products  CUR where 
the factors $C$ and $R$ are made up of $l$ columns and $k$ rows
of an input matrix, respectively, 
 the  factor $U$ has  size  $k\times l$, 
 \begin{equation}\label{eqklmn}
r\le k\le m,~r\le l\le n,~kl\ll mn,
\end{equation}
and $r$ denotes a target rank,\footnote{``In practice, the target rank is rarely known in advance. $\dots$ LRA algorithms are
usually implemented in an adaptive fashion"  -- 
with  submatrices $R$ and $C$
recursively expanded ``until the error norm satisfies the desired tolerance"
\cite[Sec. 2.4]{TYUC17}.} presumed to be low.  The $k\times l$ submatix $G$ shared by the factors $C$ and $R$ is said to be  
{\em CUR generator}. 
Given $M$ and $G$, the
factors $C$ and $R$ are immediately  available,  
and one obtains
{\em  canonical CUR} by computing the factor $U$ 
in 
about $kl$ scalars and $O(kl\min\{k,l\})$ flops
(see  (\ref{eqcnncur})).  
This computation is quasi  superfast   
 for $kl\ll mn$
 and  superfast
 for $kl\min\{k,l\}\ll mn$. 

In   
\cite{ZO18} 
Zamarashkin and
Osinsky 
proved  
 that, for any $m\times n$ matrix $M$ and any positive integer $k\le\min\{m,n\}$, there
 exists  
 a CUR LRA defined by a $k\times k$  generator and optimal up to a factor of $k+1$   under the Frobenius matrix norm.  
In \cite{CK20}  
Cortinovis and 
 Kressner  computed such a CUR  at deterministic   cost $O(kmn^3)$  for $m\ge n$
 by means of linking the task to  subset selection.
Superfast computation of a  good generator, however, is a challenge, and
  the error norms of CUR LRAs of $M$ are unbounded 
 for ill-conditioned generators $G$. 
 \\
\noindent{\bf 1.2. Our contribution.}
We extract our Part I and Part II
from  preprints \cite{PLSZa,
PLSZb}, respectively
(see further comments on this in Sec. 1.4.3), and  point out  applications
of LRA to superfast matrix norm estimation and least squares approximation of multivariate functions in Secs. \ref{snrm} and  \ref{scncl}. \\
{\bf 1.2.1. Part I: the study of CUR LRA.} 
According to    our {\bf  Main Result 1} (Thm. \ref{thm:error_estimate}),  any 
CUR LRA (no matter how  computed)
has spectral 
output error  norm $||\cdot||$  within a factor of  $(v+1)^2$
from optimal 
 for 
\begin{equation}\label{eqtheps1}v:=||U||\max\{||C||,||R||\}.\end{equation}
This  {\em universal}   upper bound implies that any
CUR LRA can fail  miserably only 
on a narrow class of inputs. In our tests  various known and novel 
superfast CUR LRA algorithms
were consistently much more accurate than this upper bound ensures and usually have spectral output error norms  within  
a factor of 30 from optima
for $n\times n$ inputs 
and $n=256,512,
1024$
(see Table \ref{tb_ranrc}).
 
Furthermore,
instead of
customary  combinations of 
C-A iterations 
 {\em Strong Rank Revealing (SRR) QR or LU factorization}
of \cite{GE96,P00},
we tested their combination
with random sampling algorithms of \cite{DMM08}, and in our tests the output accuracy of this superfast
CUR LRA was 
almost always  nearly as accurate as in the case of  the algorithms of \cite{DMM08} themselves, which run at superlinear cost and are expected to be  near-optimal
under the Forbenius norm
(see Tables  \ref{tabcadmm1} 
 and \ref{tabcadmm}).
\\
{\bf 1.2.2. Part II:  sketching LRA algorithms.}
We provide another insight into this subject  by proving a little stronger error bound  for
a large class of superfast  LRA algorithms.  First recall that with a high probability   ({\em whp})
   random  
 sketching  LRA algorithms
 of  
\cite{W14,
TYUC17} compute  near-optimal  rank-$r$ approximation\footnote
{Here and hereafter
``rank-$r$ approximation" means ``approximation of rank at most $r$".} of an 
$m\times n$ matrix $M$ under the Frobenius norm by using auxiliary  matrices
({\em sketches}) $FM$ and/or $MH$  for proper $k\times m$ and/or $n\times l$
random  {\em multipliers} 
$F$ and/or $H$, respectively,  
 $\max\{k,l\}=r+p$, $k,l,m,n,r$ satisfying 
(\ref{eqklmn}), and a
positive {\em oversampling} ({\em overestimation}) parameter $p\ll \min\{m,n\}$. 
Usually these algorithms  output LRA but not CUR LRA
(see \cite{DMM08} for an important exception), but one can
extend any LRA to CUR LRA -- at the price of increasing 
the 
spectral  norm bounds
on the output errors by  a factor of $\sqrt{1+(n-r)r)}$
(see \cite[Secs. 3.3.2 and 3.3.3]{HMT11}).

 Apart from the stage of the computation of sketches,
 random  sketching  
algorithms of \cite{HMT11}
(resp. of \cite{TYUC17})
 and their extensions in 
our Sec. \ref{sbsalg}
run superfast where $l^2\ll n$ (resp. $l^2\ll n$
and  $k^2\ll m$) and quasi superfast where $l\ll n$ (resp. $l\ll n$
and  $k\ll m$).  We compute the sketches superfast  
 as well by choosing Ultrasparse multipliers $F$ and/or $H$  (see Appendix \ref{ssrht}). This way  formal support  for near-optimality of  LRA whp is lost, but next we specify  our weaker but meaningful alternative. 

 To simplify our presentation, we restrict it to real matrices
throughout, with $M^*$ denoting transpose rather than Hermitian transpose of $M$, and to the random sketching  LRA that only uses sketches $MH$, as in 
\cite{HMT11}, while in Sec. \ref{simppre}
we extend  our results  to   algorithms that
use both sketches $FM$ and $MH$.

One can explore a variety of ways towards our goal. Maybe the simplest way is
to {\em randomly sample} $l$  columns of the matrix  $M$ or equivalently  to let $H$ be 
the  $n\times l$ leftmost submatrix 
of a random $n\times n$ permutation matrix
under a fixed randomization model;  in this case only $nl$ elements of $M$  are involved, and
the calculation itself is given "for free". We, however, 
consider the  more general class of all sketches $MH$ that can be computed superfast.
 Example \ref{exdlt}
tells us that  we cannot yield accurate LRA for worst case inputs $M$ superfast, but in our  
  numerical experiments with  Ultrasparse multipliers  $F$ and $H$ we have consistently output quite accurate  LRA for a large class
  of matrices $M$. 

Towards explaining these test  results 
and  constructing a theory,
we introduce a probabilistic structure in the  space of input
matrices $M$ rather than in the space of multipliers $H$, thus defining {\em dual 
sketching} in comparison to customary primal one, where matrix
$M$ is fixed and  multiplier $H$ is random.  

The output error norm of  sketching algorithms of \cite{HMT11} is essentially defined by the product $V_1^*H$ 
where $V_1$ is the matrix 
of the $r$ top left singular vectors of $M$, associated with its $r$ top (largest)
singular values.
Now
assume 
that $V_1^*$ is the Q factor in QR or QRP factorization of a
{\em Gaussian matrix}, 
  filled with independent identically distributed    
 normal  (standard Gaussian) random variables,
 $H$ has orthonormal columns, and the ratio $l/r$ is
 reasonably large. Then we  readily
 prove that whp the sketching algorithm outputs	LRA whose 
 both Frobenius and spectral error norms 
 are within a factor of 
 $\sqrt {1 + 16~n/l}$ from  optima; this is  essentially the first part of Thm. \ref{eqtherr} 
-- our {\bf  Main Result 2!} Its second part extends the estimate to {\em any full-rank} $H$ 
at the price of increasing the error norm bound by a factor of $
\cond_2(H):=||H||\cdot||H^+||$.
 
In Sec. \ref{serrrang} we study dual LRAs under another choice of probabilistic structure, which may be more  relevant to LRA computation in the  
real world:
  we fix a multiplier $H$
with orthonormal columns 
and let  $M:= A\Sigma B +E$ for  a perturbation matrix  $E$  and for
$A\Sigma B$ being a random pseudo  SVD of a rank-$r$ matrix of size $m\times n$, said to be  a two-sided factor-Gaussian (with expected rank 
$r$). Namely, we let $\Sigma$ be any  $r\times r$  matrix of full rank $r$ (e.g., a diagonal matrix
 with $r$ positive diagonal entries,
 as in SVD) and let $A$ and $B$ be scaled Gaussian matrices, rather than 
 matrices of singular 
vectors of SVD.\footnote{By saying  ``LRA" we assume that $r$ is much smaller than 
$\min\{m,n\}$; then we can motivate the
 above definition of random pseudo SVD by recalling
(e.g., from  
 \cite[Theorem 7.3]{E89} or 
\cite{RV09})
that $\frac{1}{\sqrt k}G$ is 
close to a matrix with  
 orthonormal columns whp for $r\ll k$ and 
an $n\times r$ Gaussian
matrix $G$.}
 Most of our study 
(in particular  Theorem \ref{eqtherr} and our Main Result 3, outlined below)
applies to a more general class of perturbed 
right factor-Gaussian matrices, of the form
$AG+E$ where $G$ is a $r\times n$ Gaussian matrix and  $A$ is an $m \times r$ matrix of full rank $r$;  transposition extends the study to perturbed 
left factor-Gaussian matrices.
 
For a matrix $M$
with a certain gap in the spectrum of its singular values\footnote{This assumption is restrictive but does not exclude  highly important kernel matrices \cite{CETW23,X24}.}  
we prove that our superfast sketching algorithms output close
rank-$r$ approximations whp on the defined probability space: namely, our bound on the output errors under this randomization model only increases
from a factor of $\sqrt {1 + 16~n/l}$
versus the optimum
 to  a factor of $\sqrt {1 + 100~n/l}$, and this is  essentially Thm. \ref{therrfctr}
-- our \textbf{Main Result 3.}
Over the same class of factor-Gaussian matrices we also strengthen the bound of
our Main Result 1 whp
(see Remark \ref{remark:factor_gaussian}).

In our tests
we sampled random multipliers $F$ and  $H$, even though
our Main Results 2 and 3 (Thms.  \ref{eqtherr} and \ref{therrfctr})
also hold where these    multipliers are fixed.
\\
{\bf 1.3. Quality of superfast dual LRA.} 
  Can we decrease the
  upper bounds of  Thms.  \ref{eqtherr} and \ref{therrfctr}
  if we  choose
 multipliers $H$ according to a proper heuristic policy?
   Such a policy can be shaped empirically,
  based on  testing various candidate multipliers $H$, but here is our semi-heuristic considerations. 
  The best-studied  LRA  
algorithms use
Gaussian   multipliers 
$F$ and $H$, 
output near-optimal LRA whp, but run at superlinear cost. 
 With sketches obtained by using {\em sparse subspace embedding} \cite{CDDRa,CFSa},
\cite[Sec. 3.3]{TYUC19},
\cite[Sec. 9]{MT20}, one devises LRA  algorithms running 
at linear computational cost.
According to \cite{L09},
 such acceleration   tends to
 make the  output accuracy  somewhat less reliable. This   problem  is partly overcome
 for {\em incoherent} matrices,   whose all entries have comparable magnitude (see \cite{CFSa} for this study and formal definition
 of incoherence). Multiplication   by the matrix
 of Subsampled Randomized Hadamard or Fourier Transform ({\em SRHT or SRFT)}  tends 
 to make a  matrix $M$ incoherent \cite{CFSa}
 but  at  superlinear cost.

Now, for our heuristic  recipe, we devise $d$-{\em Abridged SRHT multipliers}
in Appendix \ref{ssrht}.\footnote{Abridged SRFT matrices can be defined similarly to Abridged SRHT matrices and empirically have 
similar power
 (see  \cite[Sec. 13.6 and Appendix B]{PLSZ17}); unlike real Abridged SRHT matrices they involve nonreal roots of unity.} 
  For $d=\log_2(n)$ they turn into  $n\times l$ SRHT multipliers, but for small positive integers $d$  (we use  $d=3$ in our tests) they are  Ultrasparse and can  be multiplied by a dense matrix  superfast. Furthermore,
they have orthonormal columns, and so Thm.  \ref{eqtherr} 
provides dual probabilistic motivation for this otherwise heuristic choice of multipliers.
In Sec. \ref
{srndsmpl}   we provide its empirical support by testing our
  superfast extensions of  
the algorithms of 
 \cite{HMT11,
 TYUC17}.
Empirical benefits  of $d$-Abridged SRHT multipliers
were observed in \cite{PLSZb} and  apparently are still a novelty. In extensive
 numerical tests in \cite{GKLPPZa}
for $d=3$  they were at least as efficient
as  Gaussian 
 most of the time although
 caused blow up of the output errors
of LRA for a narrow class of inputs.

With this reservation we can use $d$-abridged SRHT multipliers instead of Gaussian  towards additional support of our dual model in Appendix \ref{srndoprpr}.
There we prove   that 
pre- and/or post- multiplication by Gaussian matrices turns any matrix admitting LRA   into a small-norm perturbation of a factor-Gaussian matrix. This multiplication becomes superfast if we  substitute  $3$-abridged SRHT multipliers
for Gaussian.   

At the end of Sec. 1.5 we comment on  goals  and  Abridged SRHT by-product of our numerical tests. 
 \\
\noindent {\bf 1.4. Related work.} We recall two groups of the known studies of superfast LRA.  \\
{\bf 1.4.1. Formally supported  superfast LRA.}
 Superfast algorithms
 have been  long studied 
for data sparse matrices -- 
defined by small  numbers of parameters such as ones with displacement structure
(cf. \cite{KKM79,BP94,P01,
XXG12,P15} and the bibliography therein), but  some of more recent superfast algorithms output LRAs of the families of {\em data dense matrix families},   defined by order of $mn$ parameters in the space of  $m\times n$ 
matrices. 
In this subsection
 such families 
 are specified formally,  although in some cases  in terms of parameters not available superfast.

 The   
randomized  superfast algorithms
of \cite{MW17,CETW23}
and \cite{BW18} compute LRA for  $n\times n$ Symmetric Positive
Semidefinite (SPSD) and distance matrices, respectively. The superfast algorithms of  \cite{MW17,BW18} are near-optimal whp under the Frobenius norm, while
the algorithm of \cite{CETW23}, based on random pivoted Cholesky  factorization, 
outputs oversampled LRA of size $k\times k$ for $k>r/\epsilon$ and hence is not  superfast unless  $\epsilon\gg n/r$; the algorithm, however, is highly  competitive
in practice,
particularly for SPSD kernel matrices, unlike that of \cite{MW17} 
(see \cite{AP23,
CETW23,X24}
and the references therein for
 kernel matrices and their  
applications).
 
  The
sublinear cost  
bounds  of
\cite{BW18,MW17,CETW23} do not cover 
 superlinear cost of a posteriori  estimation of  their output error norms and correctness
verification. The  superfast deterministic algorithm of \cite[Part III]{LP20}, however,  computes CUR LRA of an
$n\times n$ SPSD matrix
with both spectral and Frobenius error norms within a factor of $n$ from optimal and as by-product,  at no additional cost, estimates the output
error norm, verifying correctness.

 Chiu and Demanet in
\cite{CD13}  defined CUR generator by sampling rows and columns of an input matrix $M$ uniformly at random 
and proved that this algorithm, which we  call {\em Primitive},  is quite accurate under the Frobenius norm
whp provided that
it is directed by  a proper 
 ``regularization parameter" $\delta$ and that
  there exists an accurate LRA $XYZ\approx M$, e.g., defined by the SVD of  $M$, such that the matrices $X$ and $Z^*$ have orthonormal columns and are incoherent.

The algorithm
is  superfast
 except for the computation of proper regularization parameter $\delta$, and  a numerical example in 
 \cite[Sec. 4.1]{CD13} shows
that without this regularization, the approximation error can blow up. Likewise, in our tests in Sec. \ref{sexpr}  the output error norm of that 
algorithm performed without  regularization
  exceeded
the norms of all other tested CUR LRA algorithms by factors ranging from 100 to 5,000 in the case of $n\times n$ inputs for $n = 256, 512, 1024$. 

In \cite[Sec. 3.2]{CD13}  
  superfast  Interpolative 
Decomposition  of \cite[Sec. 3.2.3]{HMT11}
  based on 
 QR  factorization with column pivoting is added to that
 algorithm,  
 and then the authors proved that their
 estimates for
the output accuracy still hold where just one of the factors $X$ or $Z^*$ has orthonormal columns and is incoherent. 
  
 Cortinovis and 
 Ying in \cite{CY25} 
make a little   weaker assumptions
about the  matrices
$X$ and $Z^*$ -- they allow some of their column vectors  to be sparse rather than incoherent provided that  CUR LRA algorithm  incorporates another variant of the  Interpolative 
Decomposition, used earlier in  \cite{LYMHY,
X24}; we outline it
 in the next subsection. CUR LRA algorithm of
\cite{CY25}
was quite accurate in the tests in \cite{CY25} but,
 as in \cite{CD13}, is 
 not superfast -- it relies on non-canonical 
transition from CUR generator to CUR LRA based on  (\ref{eqcnn}).
  
Although the algorithms of \cite{CD13,
 CY25} are not superfast, 
analysis of their output accuracy 
is still interesting,  but   
 the provision of
incoherence
of the factors of 
an unknown LRA of  $M$
does
not hold for a large class of inputs and
  generally  cannot  be
 verified  superfast for a fixed $M$.

Our Main Results 1 and 2  (Thms. \ref{thm:error_estimate}
and \ref{eqtherr})
ensure weaker output accuracy  than   \cite{CD13,CY25}
 but 
 are  universal -- they  apply   to nearly all input matrices and either to all CUR LRAs or to all sketching LRA algorithms, respectively.
\\
{\bf 1.4.2. Heuristic superfast LRA.} 
Next we recall some  superfast LRA algorithms, whose output accuracy  was empirically stronger than    our Main Results ensure, similarly to what we observed in  
our tests. 

C-A  iterations is the most popular example of such  algorithms.
For two decades of their worldwide application they have been consistently computing accurate CUR LRA  for a large class of matrices (see  \cite{GOSTZ10,B11,OZ18,ALS24}, and the bibliography therein). 
 
 The heuristic  superfast
\cite[Alg. 2.2]{LYMHY} is    much less known but important as well. 
Its output CUR LRA is  quite accurate for a large class of  unsymmetric $n\times n$
matrices $M$
according to the test results in \cite{LYMHY}. 

 For a target rank $r$ this algorithm 
seeks 
 column and row sets that define the  factors $C$ and $R$ of CUR. Initially the two sets are empty and then are recursively increased ``a few times" based on nontrivial  extension of 
 the algorithm of \cite[Sec. 3]{CD13}.
  
Namely, one 
 first uniformly samples
$q>r$ rows of $M$ to obtain   a $q\times n$ submatrix, then applies to this submatrix
QR factorization with pivoting
as in \cite{CD13}, and finally 
 appends the first $r$ column indexes  defined by this factorization to the current set of column indexes. 
 
 One expands the current set of row indexes by applying the same algorithm to the transpose $M^*$. 
 
 When
 this recursive computation of two factors $C$ and $R$ of the CUR of $M$ is completed,   a factor $U$ is computed by means of   solving a generalized LLSP defined by $M$, $C$ and $R$. This   heuristic algorithm involves 
 order of $n^{3/2}$ memory cells
 and flops, and the authors  specify  a direction to decreasing that bound  to order of $n\log(n)$
 but admit: 
 ``We have not been able to quantify the error and success probability rigorously for this procedure
at this point." 
 
 Instead of  QR factorization
 with column pivoting,
 Jianlin Xia in \cite{X24}
 uses
SRR QR  factorization
of  \cite{GE96} and calls this technique  ``progressive
alternating direction pivoting". He combines it with his additional
techniques, in particular,
randomized
error estimation.  His algorithms, 
``numerically
shown to have nearly linear complexity" in $m+n$,  have consistently output  accurate  LRAs of  real world $m\times n$  kernel matrices $M$
in his extensive numerical
tests. He  elaborates upon
 various  techniques used for
devising his algorithms,  but not as much upon their formal support and  estimates for their output accuracy and complexity. 
He admits: ``Although a fully rigorous justification of the accuracy is lacking, we give different perspectives to motivate and
support the ideas", but it is still a challenge to   
 specify formally the large subclass of kernel matrices for which his algorithms 
 output accurate LRA. 
He cites
\cite{LP20,
PLSZ20}
as his close predecessors,
but neither his algorithms  nor \cite[Alg. 2.2]{LYMHY} significantly overlap with our current ones for CUR LRA and involve no random sketching at all. \\
{\bf 1.4.3. Our old preprints.}
In Parts I and II  we  extend the unpublished 
preprints \cite{PLSZ16,
PLSZ17,
PLSZa,
PLSZb}, which cite  \cite{PQY15,
PZ17a,PZ17b} as their predecessors, where  matrix computations under the dual probabilistic model have been formally 
studied, apparently for the first time.\footnote{Such a  model was implicit in classical initialization of the Power method with multiplication  by a coordinate vector but
has not been formally studied or even introduced.} The cited preprints
 stated  motivation  for further study of  superfast LRA as their major goal (which we share)
and seemed to be
 too unorthodox
 when they appeared, but lately,
encouraged by the
success of \cite{X24,
CY25}, we revisited
  preprints \cite{PLSZa,
PLSZb}   and
extracted from them our
 Parts I and II, respectively.
To present our work in proper historical context, we keep the results of  the numerical tests of 2019 intact and do not extend them, e.g., to the numerically stable algorithm of 
\cite{N20},
even though this
 could have strengthened their performance.
 \smallskip
 
\noindent {\bf 1.5.
 Organization of our paper.}  
We  recall some  background material on LRA, Gaussain and factor-Gaussian matrices, and CUR LRA 
   in the next two subsections
   and Sec. \ref{scurlra},
 respectively. 
 In Part I (Secs. \ref{scurlra} -- \ref{sexpr}) we study computation of CUR LRA. In Sec. \ref{sprimc}  we briefly comment on C-A iterations and  specify some more 
 rudimentary superfast CUR LRA. In Sec. \ref{scgrbckg}
  we estimate  output errors of any CUR LRA, thus proving our Main Result 1.
    In Sec. \ref{serrrnd}  we a little strengthen it in the case of a perturbed factor-Gaussian input matrix.  We devote Sec. \ref{sexpr} to our numerical experiments
for the algorithms of Part I.
In Part II (Secs. \ref{sbsalg}
 -- \ref{srndsmpl}) we study random sketching algorithms for LRA. 
 We recall and a little extend
some known ones in Sec. \ref{sbsalg} and recall their deterministic output error bounds in Sec. 
\ref{sdetrerr}. In Sec. \ref{serrranin} we prove our Main Results 2 and 3, on dual output 
 error bounds. 
  In Sec. \ref{srndsmpl} 
 we cover  numerical tests of the algorithms of Part II. 
 We devote Sec. \ref{scncl}
to conclusions.
 In Appendix  \ref{shrdin}
    we specify some 
   small families of matrices on which  any  LRA fails miserably
unless it  involves all  entries of $M$.
 In Appendix B we first comment on  some random sketching LRA and then define 
Abridged SRHT  Ultrasparse
 matrices. 
 In Appendix \ref {snrmg} we recall the estimates for the norms of a Gaussian matrix and its pseudo inverse. In Appendix \ref{srndoprpr} we prove that pre- and post-multiplication by Gaussian matrices turns any rank-$r$ matrix into a
 factor-Gaussian matrix having expected rank $r$. In Appendix \ref{sncl} we prove  necessary and sufficient
 condition for having CUR decomposition
 of a matrix.
 
 In Secs. \ref{sprimc} 
 and \ref{sbsalg}  
 we  
  describe large families   
  of LRA algorithms
covered by our Main Results. 
  In Secs. \ref{sexpr}   
 and \ref{srndsmpl}   
we   tested numerically
some specified  algorithms 
from these families, which supported  our Main Results
empirically.  This
was the only purpose of our numerical tests, but the observation of potential benefits of using $d$-Abridged SRHT multipliers has come as a by-product of that study.
\\
{\bf 1.6. Some background for LRA.}
  In this subsection,
  we write  $|||\cdot|||$    to unify the spectral norm  $||\cdot||$ and the Frobenius norm
  $||\cdot||_F$,  by following \cite[Thm. 9.1]{HMT11}.
  
Hereafter   $\mathbb R^{p\times q}$  denotes the class of $p\times q$ 
  real  matrices.  
  
For positive integers $k,l,m,n,r$ satisfying  (\ref{eqklmn}),
the triplet  
of matrices \begin{equation}\label{eqXYZkl}
X\in \mathbb R^{m\times k},~ Y\in \mathbb R^{k\times \ell},~ Z\in \mathbb R^{n\times \ell},~r\le\min\{k,\ell\},
\end{equation}
as well as 
  the product $XYZ^*$ represent
a rank-$r$ approximation of a matrix  $M\in \mathbb C^{m\times n}$
with  
  the error matrix
 $E:=M-~XYZ^*$. 
 One can  replace the triplet
$\{X,Y,Z^*\}$ with the pair $\{XY,Z^*\}$
of 
 $m\times \ell$ and $\ell\times n$ matrices or
  $\{X,YZ^*\}$ 
of  $m\times k$ and $k\times n$
matrices. 

An important 3-factor LRA of $M$ is 
its $r$-{\em top
SVD}
$M_{r}=U_{r}\Sigma_{r} V_{r}^*$
 for a diagonal matrix  $\Sigma_{r}=
 \diag(\sigma_j)_{j=1}^{r}$ of the $r$ largest singular values of $M$ and
  two  matrices $U_{r}$ and $V_{r}$ of the $r$ associated left and right singular vectors, 
  respectively.
 
  $M_{r}$ is said to be
  the $r$-{\em truncation} of $M$. One can  obtain it from a matrix $M$
 by setting to 0 all its  singular values but the $r$ top (largest) ones. $r$-top SVD of $M$ is  compact SVD of $M_r$.
 
 $M^+$ denotes the Moore--Penrose pseudo inverse of $M$. 
     
\begin{theorem} \label{thtrnc} (Eckart-Young-Mirsky Theorem, see {\rm \cite[Thm. 2.4.8]{GL13}.)} 
It holds that $$\tau_{r+1}(M):=
\min_{N:~\rank(N)=r} |||M-N|||=|||M-M_{r}|||$$ 
under both spectral and Frobenius norms:
$\tau_{r+1}(M)=\sigma_{r+1}(M)$ 
under the spectral norm, and 
$\tau_{r+1}(M)=\sigma_{F,r+1}(M):=\sqrt{\sum_{j> r}\sigma_j^2(M)}$ 
under the Frobenius norm.
\end{theorem}

 \begin{lemma}\label{lehg} {\rm (The norm of the pseudo inverse of a matrix product, see, e.g., \cite{GTZ97}.)}
Suppose that $A\in\mathbb R^{k\times r}$, 
$B\in\mathbb R^{r\times l}$
and the matrices $A$ and 
$B$ have full rank 
$r\le \min\{k,l\}$.
Then
$$|||(AB)^+|||
\le |||A^+|||\cdot|||B^+|||.$$ 
\end{lemma}


 \begin{theorem}\label{thsngr} {\rm [The impact of a perturbation of a matrix on its singular values (see \cite[Cor. 8.6.2]{GL13}).]}
 For a pair of ${m\times n}$ matrices $M$ and $M+E$ it holds that
 $$|\sigma_j(M+E)-\sigma_j(M)|\le||E||~{\rm for}~j=1,\dots,\min\{m,n\}. $$
  \end{theorem} 
  
 \begin{lemma}\label{lemma:pert_sing_space} {\rm (The impact of a perturbation of a matrix on its singular space, adapted from \cite{W72}, 
 \cite[Thm. 6.4]{S73}, \cite[Thm. 1]{GT16}.)}
Let $M$ be an $m\times n$ matrix of rank $r < \min (m, n)$ where
\begin{eqnarray*}
M = 
\begin{bmatrix}
U_r & U_{\perp}
\end{bmatrix}
~
\begin{bmatrix}
\Sigma_r & 0 \\
0 & 0
\end{bmatrix}
~
\begin{bmatrix}
V_r^T \\
V_{\perp}^T
\end{bmatrix}
\end{eqnarray*}
is its SVD, and let $E$ be a perturbation matrix  such that 
\begin{eqnarray*}
\delta = \sigma_r(M) - 2~||E||_2 > 0
~{\rm and}~ 
||E||_F \le \frac{\delta}{2}.
\end{eqnarray*}
Then there exists a matrix  such that  $P\in\mathbb R^{(n-r)\times r}$,  $||P||_F < 2~\frac{||E||_F}{\delta} < 1$, and the columns of the matrix
$\tilde{V} = V_r + V_{\perp}P$ span the right leading singular subspace 
of $\tilde M = M + E$.  
\end{lemma}

\begin{remark}\label{remark:pert_sing_space}
Matrix $\tilde{V}$  above  does not necessarily have orthogonal  columns, but one can 
 readily prove that the matrix 
$(V_r + V_{\perp}P)(I_r + P^TP)^{-1/2}$
has orthonormal columns, {\em i.e.}, $(I_r + P^TP)^{-1/2}$ normalizes $\tilde V$ (see  \cite{W72,S73,GT16}).
\end{remark}  


\noindent{\bf 1.7. Background on Gaussian and factor-Gaussian matrices.}
{\em Constant matrices} are filled with constants, unlike random matrices, filled with random variables.



\begin{theorem}\label{thrnd} {\em [Non-degeneration of a Gaussian  
matrix.]}
Suppose that 
$M\in\mathbb R^{p\times q}$ is a constant matrix, $r\le\rank(M)$,
and $F$ and $H$ are $r\times p$ and $q\times r$ independent  Gaussian 
matrices, respectively.
Then  
the matrices $F$, $H$, $FM$,  and $MH$  
have full rank $r$ 
with probability 1.
\end{theorem} 
\begin{proof}
 Rank deficiency of matrices 
$F$, $H$, $FM$, and $MH$ 
is equivalent to turning into 0
the determinants
$\det(FF^*)$, $\det(H^*H)$,
$\det(FM(FM)^*)$, and $\det((MH^*)MH)$, 
respectively. 
The claim follows because 
 these equations define algebraic varieties of  lower
dimension in the linear spaces of the entries, considered independent variables (cf., e.g., \cite[Prop. 1]{BV88}).
\end{proof}

 \begin{remark}\label{refllrnk}
Events that  occur with probability 0 are  immaterial
 for  our probability estimates, and
 hereafter we say that a  matrix  has  rank $r$   if it is has rank $r$ with probability 1.  \end{remark}
 



\begin{lemma}\label{lepr3} {\rm  [Orthogonal Invariance.] \cite[Theorem 3.2.1]{T12}}.
Suppose that  
$G$ is an $m\times n$ Gaussian 
matrix, 
$k\le \min\{m,n\}$ is a positive integer, and $S\in\mathbb R^{k\times m}$ and  $T\in\mathbb R^{n\times k}$ are constant matrices,  
 having orthonormal rows and columns, respectively.
Then $SG$ and $GT$ are random matrices having  distribution of  $k\times n$ and $m\times k$ Gaussian  random matrices, respectively.
\end{lemma}



\begin{definition}\label{deffctrg}  {\em [Factor-Gaussian matrices.]} 
Let 
$A\in \mathbb R^{m\times r}$, 
$B\in \mathbb R^{r\times n}$, and
$C\in \mathbb R^{r\times r}$ be 
three constant well-conditioned matrices of full rank $r<\min\{m,n\}$.
Let $G_1$ and $G_2$ be $m\times r$ and $r\times n$ independent Gaussian 
matrices, respectively.
Then   
$G_1B$, $AG_2$, and 
$G_1 C G_2$
are {\em left},
 {\em right}, and  
 {\em two-sided factor-Gaussian  
 matrices
 of  rank} $r$, respectively.
\end{definition} 
 
 \begin{theorem}\label{thfctrg}
 The distribution
of  two-sided $m\times n$ factor-Gaussian matrices 
$G_{m,r} C G_{r,n}$
does not change if in its definition we replace
the factor $C$ by an 
appropriate
diagonal matrix 
$\Sigma=(\sigma_j)_{j=1}^r$ such that  
$\sigma_1\ge \sigma_2\ge \dots\ge \sigma_r>0$.
\end{theorem}
\begin{proof}
Let $C=U_C\Sigma_C V_C^*$ be SVD.
Then $A=G_{m,r}U_C$ 
and $B=V_C^*G_{r,n}$
have distributions of independent $m\times r$ and $r\times n$ Gaussian matrices, respectively, 
by virtue of Lemma \ref{lepr3}. 
Hence
$G_{m,r} C G_{r,n}$ has the same distribution as 
$A\Sigma_C B$.
\end{proof} 

  
\medskip 
  
{\bf \Large PART I: COMPUTATION OF CUR LRA}

\section{CUR decomposition
and CUR LRA}\label{scurlra}

 
    
 For two sets $\mathcal I\subseteq\{1,\dots,m\}$  
and $\mathcal J\subseteq\{1,\dots,n\}$,  define
the submatrices
$$M_{\mathcal I,:}:=(m_{i,j})_{i\in \mathcal I; j=1,\dots, n},  
M_{:,\mathcal J}:=(m_{i,j})_{i=1,\dots, m;j\in \mathcal J},~{\rm and}~ 
M_{\mathcal I,\mathcal J}:=(m_{i,j})_{i\in \mathcal I;j\in \mathcal J}.$$
   
Given an $m\times n$ matrix $M$ of rank 
$r$
and its  nonsingular 
$r\times r$ submatrix
$G=M_{\mathcal I,\mathcal J}$, one can readily verify that 
$M=M'$ for
\begin{equation}\label{eqcurd} 
M'=CUR,~C=M_{:,\mathcal J},~
  U=G^{-1},~G=M_{\mathcal I,\mathcal J},~{\rm and}~R=M_{\mathcal I,:}.  
\end{equation}   
We call $G$ the {\em generator} and   
 $U$ the {\em nucleus} of {\em CUR decomposition} of $M$ where $X=C$, $Y=U$, $Z=R^*$,  and $E=O$).








CUR decomposition is extended to {\em CUR approximation}  of
a matrix $M$,  although the approximation
$M'\approx M$
 for $M'$ of (\ref{eqcurd}) tends to  be poor where the generator $G$ is 
 ill-conditioned.
 \begin{remark}\label{recgr} 
The pioneering  
 papers of 1995-2001 (cf. \cite{GTZ97,GT01} and the references in \cite{OZ18}) define CGR  approximations having 
 nuclei  $G$; ``G" can
 stand, say, for
``germ". We use the acronym CUR, which is more customary in the West. 
 ``U" can stand, say, for ``unification factor", but notice the alternatives of CNR, CCR, or CSR with    
$N$, $C$, and $S$ standing for {\em ``nucleus",
``core", and ``seed"}.
\end{remark}

 
 

By generalizing  (\ref{eqcurd}) we allow to use
$k\times l$ CUR generators
for $k$ and $l$ satisfying (\ref{eqklmn})
and to choose any $l\times k$ nucleus
$U$ for which the error matrix
$E=CUR-M$ has smaller norm.


Given two matrices $C$ and $R$,   
 the  minimal Frobenius 
 error norm  of CUR LRA 
$$||E||_F=||M-CUR||_F\le ||M-CC^+M||_F+||M-MR^+R||_F$$ is reached   for the nucleus
\begin{equation}\label{eqcnn}
U=C^+MR^+
\end{equation}
(see \cite[Eqn. (6)]{MD09}).  
We, however, cannot compute such a nucleus   superfast and instead seek   
 canonical CUR LRA (cf. \cite{DMM08,OZ18}) whose
 nucleus is the Moore-Penrose pseudo inverse of the $r$-truncation of a given CUR generator:
\begin{equation}\label{eqcnncur}
U:=G_{r}^+. 
\end{equation}
Given a generator $G$ we can compute $G_{r}^+$  by using about $kl$
scalars and $O(kl\min\{k,l\})$ flops.

We do not use the following result, but it can be of independent interest. We prove it in Appendix \ref{sncl}.
 
  \begin{theorem}\label{thncl} {\rm [A necessary and sufficient criterion for CUR decomposition.]}
   Let $M'=CUR$ be a  canonical CUR of $M$ for
 $U=G_{r}^+$, $G=M_{\mathcal I,\mathcal J}$. Then $M'=M$ if and only if $\rank(G) = \rank(M)$. 
\end{theorem}
 
  
\section{ Some superfast CUR LRA algorithms}\label{sprimc}
  

Next we present some superfast CUR LRA algorithms
as heuristic. Their crude output accuracy estimates
in Thm.  \ref{thm:error_estimate}  
 are greatly superseded in   our numerical tests of their specified versions in Sec. \ref{sexpr}.  

\subsection{Primitive and Cynical algorithms}\label{scyn}
  
        
 Given five integers $k$, $l$, $m$, $n$, and  $r$ satisfying (\ref{eqklmn}) and an $m\times n$
matrix $M$, 
sample  uniformly
 or fix  a pair  
 of  sets
 $\mathcal I$ and $\mathcal J$ of $k$ row and $l$ column indexes, respectively, to define 
 a CUR generator $G$,
 compute its $r$-truncation $G_{r}$, build on it
 a canonical CUR LRA of $M$, and  
 call this LRA algorithm 
 {\bf Primitive}.
 Given the sets 
 $\mathcal I$ and $\mathcal J$, the algorithm amounts to computing the $r$-truncation 
 $(M_{\mathcal I,\mathcal J})_{r}$
 of the matrix 
 $M_{\mathcal I,\mathcal J}$ and its pseudo inverse 
$((M_{\mathcal I,\mathcal J})_{r})^+$; this involves about $kl$ scalars and
 $O(kl\min\{k,l\})$ flops. 

The following CUR LRA algorithm
  (we call it   {\bf Cynical})\footnote{Here we allude to  the benefits of the austerity and simplicity of primitive life,   
advocated by Diogenes the Cynic, and not to shamelessness and  distrust associated with modern  cynicism.}  first fixes or samples uniformly a 
larger
submatrix $M'$ of $M$   and then compresses it into a $k\times l$ CUR generator 
of $M$ by applying one of the algorithms
of  \cite{GE96,
P00,DMM08}.

  \begin{figure}[ht] 
\centering
\begin{tikzpicture}

\draw (0,0) rectangle (4,3);

\fill[pattern=north east lines] 
    (1.2,1.25) rectangle (2,2.05);
\draw (1.2,1.25) rectangle (2,2.05);

\def\bsize{0.3}
\fill[black] 
    (1.4,1.45) rectangle  (1.65,1.7);

\end{tikzpicture}
\caption{A cynical CUR algorithm (the stripes mark a $p\times q$ submatrix; a $k\times l$ CUR generator is shown in black).
}\label{fig3} 
\end{figure}

\begin{algorithm}\label{algcnc}
 
For  an $m\times n$  matrix $M$,  a target rank $r$, and four  integers $k$, $l$, $p$, $q$ such that   
\begin{equation}\label{eqklmnpqr}  
0<r\le k\le p\le m,~r\le l\le q\le n,~{\rm and}~ kl<pq,
\end{equation}
uniformly sample  
a pair  
 of  sets
 $\mathcal I$ and $\mathcal J$ of $p$ row and $q$ column indexes, respectively,  compute a $k\times l$ CUR generator
 $G_{k,l}$
for the  $p\times q$ submatrix
$M_{\mathcal I,\mathcal J}$ 
by applying to it one of the algorithms of 
\cite{GE96,
P00,DMM08}, compute the $r$-truncation $G=G_{k,l,r}$, and 
  build on it   CUR LRA of $M$.
\end{algorithm}
  The algorithm uses
about $pq$ scalars and
$O(pq\min\{p,q\})$ 
 flops and turns into the Primitive algorithm for $p=k$ and $q=l$.
  

\subsection{Cross-Approximation (C--A) iterations}\label{scait}

We refer to \cite{GOSTZ10,
B11,OZ18, ALS24}, and the references therein for variation and performance of C-A iterations.   
 In particular their maxvol version in \cite{GOSTZ10}  
 involves $O((m+n)r)$ flops per C-A iteration seeking an $r\times r$ CUR generator for an $m\times n$
 matrix.
 Extensive effort has been invested into estimation of the overall number of C-A iterations required for convergence (see \cite{GT01,GOSTZ10,OZ18,
ALS24} and the references therein), but  adequate support of their superfast  empirical performance is still a challenge.\footnote{Goreinov et al proved in \cite{GOSTZ10}  that every C-A step of maxvol  monotone increases the volume  $|\det(G)|$,
and  reasonable bounds
on the output error norms
at a global maximum of the volume have been proved in \cite{GT01,OZ18,ALS24}.
Proving that maxvol yields global maximum superfast  is still a challenge, however, although such
bounds hold whp where an input matrix 
$M$ lies  close to a factor-Gaussian matrix  (see \cite{PLSZa}).} 
  
{\bf Initialization recipes.} 
{\em Adaptive 
C-A iterations} 
 are  initialized based on
 Gaussian elimination with complete pivoting
  (at a superlinear cost) combined with dynamic
 search for gaps in the spectrum of the singular values of $M$.
Empirically the iterations converge very fast for a large
 and important class of inputs coming from the study of ODE and PDE (cf. \cite{B11} and the references therein).    Maxvol variant of C-A iterations of \cite{GOSTZ10} 
for $q=p=k=l=r$
is initialized   by using  $O(nr^2)$ flops
and still   empirically converges fast for a large class of inputs, although the authors admit 
that their 
initialization  recipe
may not work and then suggest  using  costly Gaussian elimination with complete  pivoting. For heuristic initialization one can alternatively try using the outputs of the algorithms of
  \cite{LYMHY,
  X24,
CY25}.
 

\subsection{Maxvol
for rank-1 approximation  and norm estimation}\label{snrm}

 For a vector
${\bf v}=(v_i)^{n}_{i=1}\in\mathbb R^n$ and a
matrix $M=(m_{i,j})_{i,j=1}^{m,n}\in\mathbb R^{m\times n}$ 
   write
\begin{equation}\label{eqnrms}
||{\bf v}||_{\infty}:=
\max_{i=1}^n|v_i|,~||{\bf v}||_1:=\sum_{i=1}^n|v_i|,~
||M^*||_{\infty}:=||M||_1:=\max_{j=1}^n\sum_{i=1}^m|m_{i,j}|.
\end{equation}
Let $|M|:=\max_{i,j=1}^{m,n}|m_{i,j}|$ define 
the inconsistent matrix  norm\footnote{E.g., $|M^2|=n>|M|^2=1$ for $M=(1)_{i,j=1}^{n,n}$.}  $|M|:=\max_{i,j=1}^{m,n}|m_{i,j}|$,
 called Chebyshev's norm  in \cite{OZ18}, such that
 \begin{equation}\label{eqnrm20} 
|M|\le ||M||\le \sqrt{mn}~|M|,~|M|\le|||M||_1\le m|M|,~|M|\le||M||_{\infty}\le n|M|.
\end{equation}One can compute
the norms $|M|$, $||M||_1$,  $||M||_{\infty}$ by using  $3mn-3$ additions and comparisons.

For $p=q=r
=1$ maxvol computes an absolutely largest entry of $M=
(m_{i,j})_{i,j=1}^{m,n}$ in its fixed 
 column or row
 (in its ``vertical" 
 or ``horizontal" step), respectively, 
by using
$m-1$ or $n-1$ comparisons, respectively. 
Maxvol stops at the entry $m_{i,j}$ that maximizes 
both $|m_{g,j}|$
over all $g$ and 
$|m_{i,h}|$
over all $h$.
For a large class of matrices $M$ this is also the global maximum
$|m_{i,j}|=|M|$.

In that case 
Eqns. (\ref{eqnrms})
and (\ref{eqnrm20})  imply that the 1-norms of the $j$th column vector $C:=
(m_{g,j})_{g=1}^m$ and the $i$th row vector $R:=
(m_{i,h})_
{h=1}^n$ provide crude
 estimates for  the 1-norm and the $\infty$-norm of $M$, respectively, being off by at most factors of $m$ and $n$, respectively.

Furthermore,
write  $U:=1/m_{i,j}$,   obtain rank-1 approximation  of $M$ with the matrix $CUR$, and recall from   the result of \cite{ALS24}, extending \cite{GT01,
OZ18}, that
$$|M-CUR|\le \frac{2\sigma_1(M)\sigma_2(M)}{(\sigma_1^2(M)+\sigma^2_2(M))^{1/2}}.$$
 Thm. \ref{thtrnc}
and bound (\ref{eqnrm20})
combined imply that
 this 
approximation
 is optimal under the spectral norm up to a factor of $2\sqrt{mn}$.


\section{A posteriori errors of a canonical CUR LRA}\label{scgrbckg} 


Example \ref{exdlt} shows that {\bf a posteriori errors} of an LRA given by a matrix $M'\approx M$ cannot be estimated
 superfast  for worst case inputs, but next we prove such a bound where $M'$ is given by a triplet $\{C,U,R\}$.

\subsection{ Error Estimation: an Outline and the Statement} 
 We estimate the error norm of CUR rank-$r$ approximation of $M$ by comparing it with CUR
 decomposition of a nearby rank-$r$ matrix $M'$ defined by the same sets of row and column indexes.

\begin{outline}\label{pr1} {\rm [Error Estimation for  a Canonical CUR LRA.]} 
   \begin{enumerate}
  \item
  Consider (but do not compute) an auxiliary 
  $m\times n$ matrix $M'$ of rank $r$ that approximates the  matrix  $M$ within a fixed norm bound $\epsilon$ such that 
 \begin{equation}\label{eqmm'} 
\sigma_{r+1}(M)\le  ||M-M'||\le  \epsilon.
\end{equation} 
  [E.g.,
  we can choose $M'=M_r$ and $=\epsilon:=\sigma_{r+1}(M)$.]
\item
For the matrices $M$ and  $M'$ fix  two 
 row and column 
index sets
$\mathcal I$ and $\mathcal J$, respectively,
and define  $k\times l$
generators  $G=M_{\mathcal I,\mathcal J}$ and $G'=M_{\mathcal I,\mathcal J}'$,   nuclei $U=G_{r}^+$ and $U'=G_{r}'^+$,
and canonical  CUR approximation 
 $M\approx CUR$ and  decomposition $M'=C'U'R'$. 
  \item
  Observe that
\begin{equation}\label{eqwc} 
|| M- C U R||\le
|| M-M'||+||M'- C U R||\le \epsilon + ||C'U'R'- C U R||.
\end{equation}
\item%
Bound the norm
$||C'U'R'- C U R||$ in terms of the values
$\epsilon$, $||C||$, $||U||$, and $||R||$.
\end{enumerate}
\end{outline}  

Next we  
elaborate upon step 4 provided 
that we have already performed steps 1 -- 3. 

\begin{theorem}\label{thm:error_estimate}
Given an $m\times n$ matrix $M$, a $k\times l$ matrix $G := M_{\mathcal{I}, \mathcal{J}}$, a positive
$\epsilon\le \sigma_{r}(G)$, a positive integer $r < \min\{k, l\}$, and an $m\times n$ rank-$r$ matrix $M'$
satisfying (\ref{eqmm'}), fix $v$ of
(\ref{eqtheps1}), 
write 
\begin{equation}\label{eqv}
C:=M_{:,\mathcal{J}},~R:=M_{\mathcal{I}, :},~ 
 U = G_{r}^+,
\end{equation}

\begin{equation}\label{eqzt}
\zeta := 
\begin{cases}
\sqrt{2} & {\rm for} ~~r = \min\{k, l\} \\
(1+\sqrt{5})/2 &{\rm for}  ~~r < \min\{k, l\}
\end{cases}
\end{equation}
and assume that $\theta:=\epsilon ||U||<1$.
(We can satisfy this bound by scaling $M$  without changing $v$.) Then  it  holds that
 \begin{equation}\label{eqn:error_estimate}
    || M - CUR || \le (v+1)\Big(\frac{2\zeta}{1-\theta} ~(v+1)+ 2\Big)\epsilon,~{\rm where}~2\zeta\le 
1+\sqrt{5}. \end{equation}
\end{theorem}

 In view of (\ref{eqwc}) we only need to estimate the norm $||C'U'R' - CUR||$.

\subsection{The first bound on
$||C'U'R' - CUR||$}
\begin{lemma}\label{thwc}  
Fix  five integers  
$k$, $l$, $m$, $n$, and $r$ such that $r\le k\le m$
and $r\le l\le n$,  an  $m\times n$ matrix $M$,  
its rank-$r $ approximation $M'$  satisfying
(\ref{eqmm'}), 
and canonical CUR LRAs
$$ M\approx C U R~{\rm and}~M'=C'U'R'$$
defined by the same pair of 
 index sets
$\mathcal I$ and  $\mathcal J$ of cardinality $k$ and $l$, respectively, such that 
$$C:=M_{:,\mathcal J},~
R:=M_{\mathcal I,:},~U=G_{r }^+,~
~ C':= M'_{:,\mathcal J},~R':= M'_{\mathcal I,:},~U'=G_{r}'^+,$$
$$ G=M_{\mathcal I,\mathcal J},~{\rm and}~G'=M'_{\mathcal I,\mathcal J}.$$  
 Then   
$$||C' U' R'- C U R||\le 
(||R||+||C'||)~|| U||~\epsilon+  
||C'|| ~||R'||~||U'- U||.$$
\end{lemma}
\begin{proof}
Notice that   
$$ C U R-C'U'R'=( C-C')  U R+
C'  U( R-R')+ 
C'( U-U')R'.$$
Therefore
\begin{align*}
||  C U R-C'U'R'|| &\le
|| C-C'||~|| U||~|| R||+  
||C'||~|| U||~|| R-R'||+
||C'||~|| U-U'||~||R'||.
\end{align*}
Substitute the bound
$\max\{||C-C'||,||R-R'||\}\le ||M-M'||\le \epsilon$.
\end{proof}

\subsection{Estimation of the norm  $||U-U'||$}\label{spstr1}

 
  

\begin{lemma}\label{leprtpsdinv}   
Under the assumptions of 
Lemma \ref{thwc}
 we have
$$
||U - U'|| \le 2\zeta~||U||~||U'||~\epsilon~
  {\rm for}~\zeta~{\rm of}~(\ref{eqzt}).$$ 
\end{lemma}

\begin{proof} 
Recall that $\rank(G') = \rank(M') = \rank(G_{r}) = r$ and that $$||G' - G_{r}|| \le ||G' - G|| + ||G - G_{r}|| \le ||M' - M|| + \sigma_{r + 1}(M) \le 2\epsilon.$$ Then apply \cite[Thm. 2.2.5]{B15} for $A=G,~B=G'$. 
\end{proof}



\begin{lemma}\label{leprtinv}
Under the assumptions of Lemma \ref{leprtpsdinv}
let (\ref{eqtheps1}) hold. Then
 $||U'||\le \frac{||U||}{1-\theta}$.
\end{lemma}

\begin{proof}
Thm. \ref{thsngr} implies that  
$$\sigma_{r}(G') \ge \sigma_{r}(G)-\epsilon=(1-\epsilon/\sigma_{r}(G))\sigma_{r}(G)=(1-\epsilon~||U||)\sigma_{r}(G).$$   
Substitute (\ref{eqtheps1}) and obtain that 
$\sigma_{r}(G') \ge
(1-\theta)\sigma_{r}(G)>0$.
 Hence 
$$||U'||=\frac {1}{\sigma_{r}(G')}\le \frac{1}{1-\theta}\cdot \frac {1}{\sigma_{r}(G)} = \frac{||U||}{1 - \theta}.$$  
\end{proof}

\begin{corollary}
\label{couu'}
Under the assumptions of Lemma \ref{leprtpsdinv}
it holds that
$
||U - U'|| \le \frac{2\zeta}{1-\theta}~||U||^2~\epsilon$.
\end{corollary}

\subsection{Proof of Thm. \ref{thm:error_estimate}}

Combine bound (\ref{eqwc}) and  Lemma \ref{thwc}  and 
obtain
\begin{equation}\label{eqmcur}
||M-C U R||\le \epsilon + (||R||+||C'||)||U||\epsilon+||C'||~||R'||~||U-U'||.
\end{equation}

Recall (\ref{eqmm'}) and obtain
$||C'||\le ||C|| + \epsilon$. 

Combine  (\ref{eqtheps1})
and Cor.    
\ref{couu'}  and deduce that $$||C'||~||U-U'||\le (||C||+\epsilon) \frac{2\zeta}{1-\theta}~||U||^2\epsilon<\frac{2\zeta}{1-\theta}(||C||~|U||^2+||U||)\epsilon.$$ 
Hence
$$||C'||~||R'||~||U-U'||<\frac{2\zeta}{1-\theta}(||R'||~||C||~||U||^2+||R'||~||U||)\epsilon.$$
Substitute  relationships  (\ref{eqv}) and  obtain
\begin{equation}\label{eqcru}
||C'||~||R'||~||U-U'||<
\frac{2\zeta}{1-\theta}(v+1)^2
\epsilon.
\end{equation}

Likewise, deduce from the bounds $||C'||\le ||C|| + \epsilon$ and 
 $\epsilon||U||<1$ that 
$$(||R||+||C'||)||U||\le 
(||R||+||C||+\epsilon)||U||< (||R||+||C||)||U||+1.$$

Recall the bound on $v$ from
(\ref{eqv})
and obtain that
$(||R||+||C'||)||U||<2v+1$.

 Combine this bound with Eqns. (\ref{eqmcur})  and (\ref{eqcru})
 and  obtain the theorem.




  
\section{The errors of CUR LRA of a perturbed factor-Gaussian matrix}\label{serrrnd} 

 In this section we 
prove that   a canonical CUR LRA of a matrix $M$ 
 with a 
 generator $G = M_{\mathcal{I}, \mathcal{J}}$ 
is quite accurate whp for any fixed pair of sufficiently  large  sets 
$\mathcal{I}$ and $\mathcal{J}$ of rows and columns
indexes if $M$ is close to  a two-sided factor-Gaussian 
 matrix $M'$ of low rank.
Unlike Def. \ref{deffctrg}  we use  character $H$ for Gaussian matrices, leaving character $G$ for CUR generators.

\begin{theorem}\label{thfgcur}  
 Given two independent Gaussian matrices $H_1\in\mathbb R^{m\times \rho}$ and   $H_2\in\mathbb R^{\rho\times n}$ and  a $\rho\times \rho$ (well-conditioned)
matrix $\Sigma$ of full rank   
 $\rho$
 such that
 \begin{equation}\label{eqn:fgpreq}
    k < m,~l < n,~\textrm{and}~ \min\{k, l\}\ge \rho + 2 \ge 4, 
\end{equation}
 write
 \begin{equation}\label{eqn:fgdfn}
    M':= H_1 \Sigma H_2,
\end{equation}
$$C':=M'_{:,\mathcal J},~ 
R':=M'_{\mathcal I,:},~
G'= M'_{\mathcal I,\mathcal J}\in \mathbb R^{k\times l}, ~{\rm and}~U':=G'^+$$ 
for fixed
 row and column 
index sets $\mathcal I$ and $\mathcal J$.
Then $C'$,  $R'$, and $G'$ are factor-Gaussian matrices such that
\begin{eqnarray*}
    \mathbb{E} ||C'|| \le (\sqrt{m} + \sqrt{\rho})(\sqrt{\rho} + \sqrt{l})\sigma_1(\Sigma)\le \beta\sigma_1(\Sigma),~~~~~~~~~ \\
    \mathbb{E} ||R'|| \le (\sqrt{n} + \sqrt{\rho})(\sqrt{\rho} + \sqrt{k})\sigma_1(\Sigma)\le \beta\sigma_1(\Sigma),~~~~~~~~~\\
    \mathbb{E}||U'|| 
   \le \frac{e^2\sqrt{kl}}{(k-\rho)(l-\rho)\sigma_\rho(\Sigma)}\le \frac{\alpha}{\sigma_\rho( 
    \Sigma)}  ~~~~~~~~~~~~~~~~~~~~~
    \end{eqnarray*}
~{\rm for}~ $e:=2.7182818
\dots$,  \begin{equation}\label{eqn:alpbeta}
    \alpha := \frac{e^2\sqrt{kl}}{(k-\rho)(l-\rho)}
    ~\textrm{and}~ 
    \beta := \max\big\{ (\sqrt{m} + \sqrt{\rho})(\sqrt{\rho} + \sqrt{l}), (\sqrt{n} + \sqrt{\rho})(\sqrt{\rho} + \sqrt{k}) \big\}.
\end{equation}

\end{theorem}
\begin{proof}
We  estimate
$\mathbb{E}||C'||$ 
and  one can 
similarly  estimate
$\mathbb{E}||R'||$.
 
 Recall that $C' = M'_{:, \mathcal{J}} =  H_1\Sigma{H_2}_{:, \mathcal{J}}$ and
  $H_1$ and $H_2$ are independent
 of one another
 and deduce that
\begin{align*}
    \mathbb{E} ||C'|| &\le \mathbb{E} \big( ||H_1||~||\Sigma||~||{H_2}_{:,\mathcal{J}}||\big)\\
    &= \sigma_1(\Sigma) \mathbb{E}||H_1|| ~ \mathbb{E}||{H_2}_{:,\mathcal{J}}||.
\end{align*}
Recall from Thm. \ref{thsignorm} that $\mathbb{E}(\nu_{p, q}) \le \sqrt{p} + \sqrt{q}$
and obtain
\begin{equation*}
    \mathbb{E} ||C'|| \le (\sqrt{m} + \sqrt{\rho})(\sqrt{\rho} + \sqrt{l})\sigma_1(\Sigma).
\end{equation*}
Next 
estimate $||U'||$. Apply Lemma \ref{lehg},
recall from Thm. \ref{thsiguna} that 
$$\mathbb{E}(\nu^+_{p, q}) \le \frac{e\sqrt{p}}{p - q},~{\rm for}~p \ge q+2 \ge 2,$$  and deduce that 
\begin{align*}
\mathbb{E}||U'|| &= \mathbb{E}||({H_1}_{\mathcal{I}, :}\Sigma{H_2}_{:, \mathcal{J}})^+||\\
&\le \mathbb{E}||{H_1}_{\mathcal{I}, :}^+||~\mathbb{E}||{H_2}_{:, \mathcal{J}}^+||~/\sigma_{\rho}(\Sigma)   \\
&\le \frac{e^2\sqrt{kl}}{(k-\rho)(l-\rho)\sigma_\rho(\Sigma)}= \frac{\alpha}{\sigma_\rho(\Sigma)}.
\end{align*}

\end{proof}
 
 
Now  assume that an  input matrix $M$ is a perturbed two-sided factor-Gaussian matrix and that two  
fixed sets $\mathcal{I}$ and $\mathcal{J}$ of its row and column indexes are sufficiently large and then
estimate
the norms  $||U||$,  $||C||$ and $||R||$.

\begin{lemma}\label{lemma:factor_gaussian_prob}
Under the assumptions of Thm. \ref{thfgcur}, let $\epsilon \le \frac{\sigma_\rho(\Sigma)}{30\alpha}$ be a small positive number and 
let $E\in \mathbb{R}^{m\times n}$ be a perturbation matrix such that $||E|| \le \epsilon$. 
Write $$M:= M' + E,~C := M_{:, \mathcal{J}},~R:= M_{\mathcal{I}, ;},~{\rm and}~U: = (M_{\mathcal{I}, \mathcal{J}})_\rho^+.$$ Then with a probability no less than 0.7 we have 
$$ {\rm (i)}~~~~~~~~~~~~~~~~~~~~~~~~~~~~~~~~~~~~~~~~~||U|| \le \frac{15\alpha}{\sigma_{\rho}(\Sigma)},~~ \epsilon ||U|| \le 1/2, ~~~~~~~~~~~~~~~~~~~~~~~~~~~~~~~$$

$$ {\rm (ii)}~~~~~~~~~~~~~~~~~~~~~~~~~~~
\max\{ ||C||, ||R|| \}||U|| \le 150\alpha\beta\sigma_1(\Sigma)/\sigma_\rho(\Sigma) + 1/2.~~~~~~~~~~~~~~~~~~$$
\end{lemma}

\begin{proof}
(i) By combining the bound of Thm. \ref{thfgcur} on the expected norm
$\mathbb  E||U'||$ with Markov's inequality,  deduce that 
\begin{equation}\label{equ'alph}
 ||U'|| \le \frac{10\alpha}{\sigma_\rho(\Sigma)} \textrm{ or equivalently } \sigma_{\rho}(M'_{\mathcal{I}, \mathcal{J}}) \ge \frac{\sigma_\rho(\Sigma)}{10\alpha}
\end{equation}
with 
a probability no less than 0.9.

Similarly to the argument of Lemma \ref{leprtinv},  deduce that $$\sigma_{\rho}(M_{\mathcal{I}, \mathcal{J}}) \ge \sigma_{\rho}(M'_{\mathcal{I}, \mathcal{J}}) - \epsilon \ge \frac{\sigma_\rho(\Sigma)}{15\alpha},$$ and then claim (i) 
follows.\\
(ii) By combining the bounds of Thm. \ref{thfgcur} on the expected norms  
$\mathbb||C'||$ and $\mathbb||R'||$ with
 Markov's inequality and the  union probability bound for random variables, 
 obtain that 
$$
\max\big\{||C||,||R||\big\} \le 10\beta\sigma_{1}(\Sigma) + \epsilon
$$
with 
  a
 probability no less than 0.8.
Then combine claim (i) with the union bound
and obtain that
\begin{equation*}
    \max\{ ||C||, ||R|| \}||U|| \le 150\alpha\beta\sigma_1(\Sigma)/\sigma_\rho(\Sigma) + 1/2  
\end{equation*}
 with 
a probability no less than 0.7.
\end{proof}

Now choose 
 any generator of CUR LRA of a perturbed factor-Gaussian matrix of rank
 $\rho$ and prove a  probabilistic upper bound on    the error of a 
resulting 
CUR LRA of 
that matrix. 

\begin{theorem}\label{thm:pfg_error_bound}
If the 
assumptions
of Lemma \ref{lemma:factor_gaussian_prob}
hold for  a perturbed factor-Gaussian matrix $M$, then 
bound (\ref{eqn:error_estimate})  on 
the approximation error norm $||M - CUR||$ for $\theta\le 1/2$, that is, 
 \begin{equation}\label{eqn:factor_gauss_error_est} 
  ||M - CUR|| \le (     4\zeta v+ 2)(v+1) \epsilon,~{\rm for}~4\zeta\le 2+2\sqrt 5,~v = 150\alpha\beta\sigma_1(\Sigma)/\sigma_\rho(\Sigma) + 1/2,
\end{equation}
holds with a probability no less than 0.6.
\end{theorem}

\begin{proof}
Combine  bounds 
$\epsilon \le \frac{\sigma_\rho(\Sigma)}{30\alpha}$
and (\ref{equ'alph}) of  Lemma \ref{lemma:factor_gaussian_prob}  
and obtain that
 $\sigma_\rho(G)\ge 3\epsilon > \epsilon$
  and therefore  (\ref{eqn:error_estimate}) hold with a probability $p_0\ge 0.9$. 
 Replace $\theta$ and $v$ with the upper bounds 1/2 and $ 150\alpha\beta\sigma_1(\Sigma)/\sigma_\rho(\Sigma) + 1/2$,  respectively,  holding with a 
 probability $p_1\ge 0.7$.
 Combine this bound with $p_0\ge 0.9$.
\end{proof}

\begin{remark}\label{remark:factor_gaussian}
Recall that 
$$\alpha := \frac{e^2\sqrt{kl}}{(k-\rho)(l-\rho)}~{\rm and}~\beta := \max\big\{ (\sqrt{m} + \sqrt{\rho})(\sqrt{\rho} + \sqrt{l}), (\sqrt{n} + \sqrt{\rho})(\sqrt{\rho} + \sqrt{k}) \big\}.$$ Furthermore, suppose that
 $m \gg k \gg \rho$, $n \gg l \gg \rho$, and 
$$    v = 150\alpha\beta\sigma_1(\Sigma)/\sigma_\rho(\Sigma) + 1/2 = O\Big(\max\Big\{\frac{\sqrt{m}}{\sqrt{k}}, \frac{\sqrt{n}}{\sqrt{l}} \Big\}\Big).$$
Now drop the smaller terms of (\ref{eqn:factor_gauss_error_est}), and obtain that its dominant part  is
\begin{equation*}
    O\Big(\max\Big\{\frac{m}{k},\frac{n}{l} \Big\}\cdot\epsilon \Big).
\end{equation*}
\end{remark}

\section{Numerical experiments 
}\label{sexpr}


\subsection{Test overview}

 
   For  both synthetic and real world input matrices
 the test
  results
  in our Table \ref{tb_ranrc} 
  show  good accuracy of CUR LRA output by just five C-A loops (with uniform initialization). Accuracy was still quite good   where we applied the  Cynical
 algorithm but significantly worse  where we applied the Primitive algorithm. 
Tables \ref{tabcadmm1} 
 and \ref{tabcadmm} 
compare the output accuracy of
the random sampling algorithms 
of \cite{DMM08}
and their superfast combination with C-A iterations.

We have performed
 the tests   
 in the Graduate Center of the City University of New York 
by applying  MATLAB. We generated
synthetic inputs by using  Gaussian matrices, which we created by using standard normal distribution function ``randn()". We	 calculated 
 $\epsilon$-ranks of   matrices for 
 $\epsilon=10^{-6}$
 by using the MATLAB's function 
 "rank(-,1e-6)",
 which only counts singular values greater than $10^{-6}$.  
       Some numerical experiments were executed  
      with software custom programmed in $C^{++}$ and compiled with LAPACK version 3.6.0 libraries. 
  
 Our tables display the  mean value of the relative 
 output errors   $ \frac{\|M - \tilde M \|}{\|M \|}$ for an input matrix $M$ and output $\tilde M:=CUR$
 and the standard deviation (std)  over 1000 runs for every class of inputs.
  

\subsection{Input matrices for LRA}\label{ststmtrcs0}

   
 We used the following two classes of 
 input matrices $M$ for testing LRA algorithms. 

\medskip

{\bf Class I} (Synthetic inputs):  
  $n\times n$
random matrices $M$ having expected 
 rank $r$, that is, 
$$M = G_1 * G_2 + 10^{-10} G_3,$$
for three Gaussian matrices $G_1$
of size $n\times r$, $G_2$ of size $r\times n$,
and  $G_3$
 of size $n\times n$. 
\medskip

{\bf Class II}:  The dense  matrices with smaller ratios of ``$\epsilon$-rank/$n$"  from the built-in test problems in Regularization 
Tools, which came from discretization (based on Galerkin or quadrature methods) of the Fredholm  Integral Equations of the first kind,\footnote{See 
 http://www.math.sjsu.edu/singular/matrices and 
  http://www2.imm.dtu.dk/$\sim$pch/Regutools 
  
For more details see Chapter 4 of the Regularization Tools Manual at \\
  http://www.imm.dtu.dk/$\sim$pcha/Regutools/RTv4manual.pdf } namely,
to the following six input families  from the Database:
 
\medskip

{\em baart:}       Fredholm Integral Equation of the first kind,

{\em shaw:}        one-dimensional image restoration model,
 
{\em gravity:}     1-D gravity surveying model problem,
 




{\em wing}:        problem with a discontinuous
 solution,

{\em foxgood:}     severely ill-posed problem, 
 
{\em inverse Laplace:}   inverse Laplace transformation.
    
\medskip 


\subsection{Four test algorithms}


 We tested  the following four algorithms
for computing CUR LRA of input matrices $M$ having
 $\epsilon$-rank $r$:
\begin{itemize}
\item
{\bf Tests 1 (The
Primitive
 algorithm for $k=l=r$):}
  Uniformly sample two index sets $\mathcal{I}$ and $\mathcal{J}$, both of cardinality $r$, then compute a nucleus 
$U=M_{\mathcal{I}, \mathcal{J}}^{-1}$ and define  CUR LRA
\begin{equation}\label{eqtsts1}
\tilde M:=CUR=M_{:, \mathcal{J}} \cdot  M_{\mathcal{I}, \mathcal{J}}^{-1} \cdot M_{\mathcal{I},\cdot}.
\end{equation} 
%
\item    
{\bf Tests 2 (Five loops of  C--A):}
Uniformly sample an  initial row index set 
$\mathcal{I}_0$ of cardinality $r$, then perform five loops of C--A
 (cf. Sec. \ref{scait}) incorporating Alg. 1 of \cite{P00}
 as a subalgorithm  that
  produces $r\times r$ CUR generators.
  At the end compute a nucleus $U$ and define 
CUR LRA  as
in Tests 1.
\end{itemize} 
 
\begin{itemize}
\item
{\bf Tests 3 (A Cynical algorithm for $p=q=4r$ and $k=l=r$):}
Uniformly sample a row index set $\mathcal{K}$ and a column index set $\mathcal{L}$,
both of cardinality $4r$. Then apply Algs. 1 and 2 of \cite{P00} 
to compute an $r\times r$ submatrix $M_{\mathcal{I}, \mathcal{J}}$ of $M_{\mathcal{K}, \mathcal{L}}$. Compute a nucleus and obtain  CUR LRA of (\ref{eqtsts1}).

\item
{\bf Tests 4 (Combination of a single 
 C--A loop with Tests 3):}
Uniformly sample a  column index set 
$\mathcal{L}$
of cardinality $4r$; then
perform a single C--A
loop: first define an index set $\mathcal{K}'$ of cardinality $4r$ and the submatrix  $M_{\mathcal{K}', \mathcal{L}}$ in  $M_{:, \mathcal{L}}$
 by applying Alg. 1 of \cite{P00}; then apply this algorithm to the matrix  $M_{\mathcal{K}',:}$ to find an index set  
$\mathcal{L}'$ of cardinality $4r$ and to define  submatrix  $M_{\mathcal{K}', \mathcal{L}'}$  in $M_{\mathcal{K}',:}$. Then proceed as in Tests 3 -- find an $r\times r$ CUR generator  $M_{\mathcal{I}, \mathcal{J}}$ in $M_{\mathcal{K}', \mathcal{L}'}$ by applying Algs. 1 and 2 of \cite{P00} and finally  compute a nucleus and  CUR LRA of $M$ based on (\ref{eqcnncur})
for $G:=M_{\mathcal{I}, \mathcal{J}}$.
\end{itemize}


\subsection{CUR LRA of the  matrices of class I}\label{sinteq0}

Table \ref{tb_ranrc} shows the summary statistics   for each of the four tests matrices of class I for $n =256, 512, 1024$ and $r = 8, 16, 32$.  In the column SVD we display the data
for the best rank-$r$ approximation given by the $(r+1)$-st largest singular value to establish the baseline for the performance of our tests.  The results of Tests 1 fall in the range $[10^{-8},10^{-7}]$,
those of Tests 2, 3 and 4
 in $[10^{-11},10^{-10}]$  in the decreasing order of accuracy. 
On the average the error bounds of Tests 2, 4, 3, and 1 exceed  the baseline
bounds in the column SVD by roughly factors of 10,
15, 30, and  over 10,000, respectively,  but even  the crudest LRAs output in  Test 1  could be used in some applications.

\setlength\tabcolsep{4 pt}
\renewcommand{\arraystretch}{1.2}
\begin{table}[h!]
\small
\begin{center}
\begin{tabular}{|c c|c c|c c|c c|c c|c c|}
\hline     
 &  &\multicolumn{2}{c|}{\bf SVD} &\multicolumn{2}{c|}{\bf Tests 1} & \multicolumn{2}{c|}{\bf Tests 2}&\multicolumn{2}{c|}{\bf Tests 3} & \multicolumn{2}{c|}{\bf Tests 4}\\ \hline
{\bf n} & {\bf r} & {\bf mean} & {\bf std}  &{\bf mean} & {\bf std}  & {\bf mean} & {\bf std}& {\bf mean} & {\bf std} & {\bf mean} & {\bf std}\\ \hline
256 & 8     & 1.01e-11 & 3.92e-13 & 1.60e-08 & 8.65e-08 & 5.94e-11 & 8.06e-12 & 1.13e-10 & 2.36e-11 & 8.23e-11 & 1.64e-11 \\ \hline
256 & 16   & 9.12e-12 & 3.09e-13 & 2.44e-07 & 5.85e-06 & 7.31e-11 & 1.02e-11 & 1.12e-10 & 2.02e-11 & 9.45e-11 & 1.51e-11 \\ \hline
256 & 32   & 7.80e-12 & 2.32e-13 & 4.82e-08 & 2.50e-07 & 8.93e-11 & 1.10e-11 & 1.13e-10 & 1.69e-11 & 1.04e-10 & 1.50e-11 \\ \hline
512 & 8     & 7.64e-12 & 2.23e-13 & 3.50e-08 & 2.72e-07 & 5.71e-11 & 7.08e-12 & 1.21e-10 & 2.67e-11 & 8.34e-11 & 1.62e-11 \\ \hline
512 & 16   & 7.06e-12 & 1.68e-13 & 1.18e-07 & 2.53e-06 & 7.08e-11 & 8.96e-12 & 1.26e-10 & 2.21e-11 & 9.98e-11 & 1.59e-11 \\ \hline
512 & 32   & 6.36e-12 & 1.37e-13 & 7.43e-08 & 7.47e-07 & 9.25e-11 & 1.14e-11 & 1.34e-10 & 1.92e-11 & 1.20e-10 & 1.73e-11 \\ \hline
1024 & 8   & 5.63e-12 & 1.13e-13 & 2.42e-08 & 2.46e-07 & 5.39e-11 & 5.89e-12 & 1.28e-10 & 2.83e-11 & 8.10e-11 & 1.55e-11 \\ \hline
1024 & 16 & 5.34e-12 & 9.23e-14 & 6.12e-08 & 6.91e-07 & 6.94e-11 & 7.68e-12 & 1.37e-10 & 2.35e-11 & 1.04e-10 & 1.73e-11 \\ \hline
1024 & 32 & 4.95e-12 & 7.55e-14 & 6.20e-07 & 1.36e-05 & 9.17e-11 & 1.06e-11 & 1.51e-10 & 2.09e-11 & 1.29e-10 & 1.87e-11 \\ \hline
\end{tabular}
\caption{Errors of CUR LRA of random matrices of class I}
\label{tb_ranrc}
\end{center}
\end{table}
 \renewcommand{\arraystretch}{1}
 \setlength\tabcolsep{6pt}

\subsection{LRA by means of random sampling
and C-A acceleration}
\label{ststc-alvrg}


 Tables  
 \ref{tabcadmm1} 
 and    
 \ref{tabcadmm}   
  display the  relative errors
  $ \frac{\|M - \tilde M \|}{\|M \|}$ of the LRA $\tilde M$  of  the matrices $M$ 
computed  in two ways:   by means of
  \cite[Alg. 2]{DMM08} (see the lines marked ``CUR"), which are expected to be near-optimal under the Frobenius morm whp,
 or eight C-A iterations combined with \cite[Alg. 1]{DMM08}  
 (see the lines marked ``C-A"). The  algorithms run at superlinear 
  cost in the former case and  are superfast in the latter case.
   The columns of the tables marked with "$\epsilon$-rank" display  $\epsilon$-rank of an input matrix.
 The columns of the tables  marked with "$k=l$" show the number of rows and  columns in a square matrix of CUR generator. 

 As we pointed out in Sec. 1.2.1, in  almost all our tests this dramatic acceleration was 
 achieved  at the price of only minor deterioration of output accuracy. 
 \medskip


\begin{table}[ht]
\centering
\begin{tabular}{|c|c|c|c|c|c|c|c|}
\hline
input & algorithm & m & n & $\epsilon$-rank & k=l & mean & std \\ \hline
finite diff & C-A & 608 & 1200 & 94 & 376  & 6.74e-05 & 2.16e-05 \\ \hline
finite diff & CUR & 608 & 1200 & 94 & 376  & 6.68e-05 & 2.27e-05 \\ \hline
finite diff & C-A & 608 & 1200 & 94 & 188  & 1.42e-02 & 6.03e-02 \\ \hline
finite diff & CUR & 608 & 1200 & 94 & 188  & 1.95e-03 & 5.07e-03 \\ \hline

baart & C-A & 1000 & 1000 & 6 & 24 & 2.17e-03 & 6.46e-04 \\ \hline
baart & CUR & 1000 & 1000 & 6 & 24 & 1.98e-03 & 5.88e-04 \\ \hline
baart & C-A & 1000 & 1000 & 6 & 12 & 2.05e-03 & 1.71e-03 \\ \hline
baart & CUR & 1000 & 1000 & 6 & 12 & 1.26e-03 & 8.31e-04 \\ \hline
baart & C-A & 1000 & 1000 & 6 & 6  & 6.69e-05 & 2.72e-04 \\ \hline
baart & CUR & 1000 & 1000 & 6 & 6  & 9.33e-06 & 1.85e-05 \\ \hline

shaw & C-A & 1000 & 1000 & 12 & 48 & 7.16e-05 & 5.42e-05 \\ \hline
shaw & CUR & 1000 & 1000 & 12 & 48 & 5.73e-05 & 2.09e-05 \\ \hline
shaw & C-A & 1000 & 1000 & 12 & 24 & 6.11e-04 & 7.29e-04 \\ \hline
shaw & CUR & 1000 & 1000 & 12 & 24 & 2.62e-04 & 3.21e-04 \\ \hline
shaw & C-A & 1000 & 1000 & 12 & 12 & 6.13e-03 & 3.72e-02 \\ \hline
shaw & CUR & 1000 & 1000 & 12 & 12 & 2.22e-04 & 3.96e-04 \\ \hline
\end{tabular}
\caption{LRA errors of Cross-Approximation (C-A)   incorporating 
\cite[Algorithm 1]{DMM08} in comparison to stand-alone CUR of 
\cite[Algorithm 2]{DMM08}.}
\label{tabcadmm1}
\end{table} 

\begin{table}[ht]
\centering
\begin{tabular}{|c|c|c|c|c|c|c|}
\hline
input & algorithm & m = n & $\epsilon$-rank & $k=l$ & mean & std \\ \hline
foxgood & C-A & 1000  & 10 & 40 & 3.05e-04 & 2.21e-04 \\\hline
foxgood & CUR & 1000  & 10 & 40 & 2.39e-04 & 1.92e-04 \\ \hline
foxgood & C-A & 1000  & 10 & 20 & 1.11e-02 & 4.28e-02 \\ \hline
foxgood & CUR & 1000 & 10 & 20 & 1.87e-04 & 4.62e-04 \\ \hline
wing & C-A & 1000  & 4 & 16 & 3.51e-04 & 7.76e-04 \\ \hline
wing & CUR & 1000  & 4 & 16 & 2.47e-04 & 6.12e-04 \\ \hline
wing & C-A & 1000  & 4 & 8  & 8.17e-04 & 1.82e-03 \\ \hline
wing & CUR &  1000 & 4 & 8  & 2.43e-04 & 6.94e-04 \\ \hline
wing & C-A &  1000 & 4 & 4  & 5.81e-05 & 1.28e-04 \\ \hline
wing & CUR &  1000 & 4 & 4  & 1.48e-05 & 1.40e-05 \\ \hline

gravity & C-A &  1000 & 25 & 100 & 1.14e-04 & 3.68e-05 \\ \hline
gravity & CUR &  1000 & 25 & 100 & 1.41e-04 & 4.07e-05 \\ \hline
gravity & C-A &  1000 & 25 & 50  & 7.86e-04 & 4.97e-03 \\ \hline
gravity & CUR &  1000 & 25 & 50  & 2.22e-04 & 1.28e-04 \\ \hline
inverse Laplace & C-A &  1000 & 25 & 100 & 4.15e-04 & 1.91e-03 \\ \hline
inverse Laplace & CUR &  1000 & 25 & 100 & 5.54e-05 & 2.68e-05 \\ \hline
inverse Laplace & C-A &  1000 & 25 & 50  & 3.67e-01 & 2.67e+00 \\ \hline
inverse Laplace & CUR &  1000 & 25 & 50 &  2.35e-02 & 1.71e-01 \\ \hline
\end{tabular}
\caption{LRA errors of Cross-Approximation (C-A)   incorporating 
\cite[Algorithm 1]{DMM08}  in comparison to stand-alone CUR of 
\cite[Algorithm 2]{DMM08}.} \label{tabcadmm}
\end{table}


 
  







\medskip 
  
{\bf \Large PART II: LRA VIA RANDOM SKETCHING}



\section{{\em LRA} by means of sketching: four  algorithms}\label{sbsalg}
 

Next we slightly generalize the sketching  {\em LRA} algorithms of \cite{HMT11,
TYUC17}
(see Remarks \ref{realg10}
 -- \ref{resmpl});
 Thms.
 \ref{thst3} and \ref{eqtherr} 
combined  formally support  output accuracy of dual version of all of them whp;
 we applied
 our tests to the algorithms of
 \cite{HMT11,
 TYUC17} accelerated by means of using 
 Abridged SRHT multipliers $F$ and $H$.
   

\begin{algorithm}\label{alg1} 
   
 
\begin{description}


\item[{\sc Input:}] 
An $m\times n$ matrix  $M$ and a
 target  rank  $r$.  


\item[{\sc Output:}] 
Two matrices $X\in \mathbb R^{m\times l}$ and 
$Y\in \mathbb R^{l\times n}$ for $r\le l\le n$ defining
an {\em LRA} $\tilde M=XY$ of $M$.


\item[{\sc Initialization:}] 
 Fix an integer $l=r+p\le n$, for $p\ge 0$, and
 an $n\times l$  matrix $H$ of full rank $l$.

 
\item[{\sc Computations:}]

\begin{enumerate}
\item 
Compute the  $m\times l$ matrix $MH$.
\item 
Fix a nonsingular $l\times l$ matrix $T^{-1}$  and output  the $m\times l$ matrix $X:=MHT^{-1}$.
\item 
Output an $l\times n$ matrix 
$Y:= {\rm argmin}_V ~|||XV-M|||=X^+M$. 
\end{enumerate}


\end{description}

 
\end{algorithm}


\begin{algorithm}\label{alg0} 
   
 
\begin{description}


\item[{\sc Input:}] 
As in Alg. \ref{alg1}.



\item[{\sc Output:}] 
Two matrices $X\in \mathbb R^{k\times n}$ and 
$Y\in \mathbb R^{m\times k}$ defining
an {\em LRA} $\tilde M=YX$ of $M$.


\item[{\sc Initialization:}] 
 Fix an integer $k=r+p\le m$, 
 for $p\ge 0,$ and
 a $k\times m$  matrix $F$ of full numerical rank $k$.

 
\item[{\sc Computations:}]

\begin{enumerate}
\item 
Compute the  $k\times m$ matrix $FM$.
\item 
Fix a nonsingular $k\times k$ matrix $S^{-1}$; then output
 $k\times n$ matrix $X:=S^{-1}FM$.
\item 
Output an  $m\times k$ matrix 
$Y:= {\rm argmin}_V ~|||VX-M|||=MX^+$.
\end{enumerate}


\end{description}

 
\end{algorithm}


\noindent 


The following
algorithm combines row and column sketching. 
  
\begin{algorithm}\label{alg01} 

 
\begin{description}


\item[{\sc Input:}] 
As in Alg. \ref{alg1}.


\item[{\sc Output:}] 
Two matrices $X\in \mathbb R^{m\times k}$ and 
$Y\in \mathbb R^{k\times n}$ defining
an {\em LRA} $\tilde M=XY$ of $M$.


\item[{\sc Initialization:}] 
 Fix two
 integers $k$ and $l$, $r\le k\le m$
 and $r\le l\le n$; fix
 two  matrices $F\in\mathbb R^{k\times m}$ 
and $H\in\mathbb R^{n\times l}$  
 of full numerical ranks
 and two nonsingular matrices $S\in\mathbb R^{k\times k}$ and $T\in\mathbb R^{l\times l}$.

 
\item[{\sc Computations:}]
1. Output the matrix  
$X=MHT^{-1 }\in\mathbb R^{m\times l}$.

 2. Compute  the matrices  
 $U:=S^{-1}FM\in\mathbb R^{k\times n}$ and
$W:=S^{-1}FX\in\mathbb R^{m\times l}$.

3. Output the $l\times n$ matrix 
$Y:={\rm argmin}_V |||WV-U|||=W^+U$.


\end{description}

 
\end{algorithm}
   


\begin{remark}\label{realg10} 
 We can obtain Algs. \ref{alg0}
  by applying Alg. \ref{alg1}
  to the transpose $M^*$.
Likewise,
by applying
Alg. \ref{alg01}
 to $M^*$ we obtain
 {\bf Alg.  6.4}.
We only study  Algs. \ref{alg1} and \ref{alg01}, but can  readily extend that study to 
Algs. \ref{alg0}
and 6.4.  \end{remark}

\begin{remark}\label{re0} 
Fix $k=r$,  $l>k$, a random  $n\times l$
matrix $H$ (e.g.,  Gaussian  or SRHT matrix), identity matrix $S$, and $T$ equal to the $R$ factor in the $QR$ factorization of $MH$. Then the matrix $X$ has  orthonormal
columns and 
Algs. \ref{alg1}
turns into Proto-algorithm of \cite[Sec. 9]{HMT11}, while Alg.
 \ref{alg01}
turns into the {\rm Generalized   
Nystr{\"o}m algorithm} 
of \cite[Eqn. (3)]{N20}.
\cite{N20} stabilizes the latter  algorithm numerically
 -- essentially by means of setting to 0 all
singular
values of the matrix $W$ exceeded by a fixed 
 $\epsilon$,  ``a modest multiple
of the unit roundoff $u$ times $||W||$". 
Our  study can be readily extended  to such a stabilized
{\em LRA} because the stabilization 
little affects  complexity of {\em LRA} and  
only improves its output  accuracy.
We can obtain
various other modifications 
of Algs. \ref{alg1} and \ref{alg01}
 by fixing other
sketching matrices $F$ and $H$.

 \end{remark}


\begin{remark}\label{resmpl}Column (resp. row) sketching turns into  column (resp. row) {\em subset selection} where $H$ in Alg. \ref{alg1}  (resp. $F$ in Alg. \ref{alg0}) is a {\em sampling matrix}, that is, a  full rank submatrix of a permutation matrix.  \end{remark}

         
   



   
\section{Deterministic output error bounds for  sketching algorithms}\label{sdetrerr}
  
  
\subsection{Deterministic error bounds of Range Finder}\label{serrrng}

                                                                                                                                                                                                                                                                                                                                                                                  Next we recall some known estimates for the errors of Alg. \ref{alg1}, to be used in the next section.
                                                                                                                                                                                                                                                                                                                                                                                  
                                                                                                                                                                                                                                                                                                                                                                                  \begin{theorem}\label{thpert1} {\rm \cite[Thm. 9.1]{HMT11}.}
Suppose that Alg. \ref{alg1}
has been applied to a matrix $M$ and let  

\begin{eqnarray*}
 M=\begin{pmatrix} U_1& U_2\end{pmatrix}
\begin{pmatrix}\Sigma_1& \\
 &\Sigma_2 \end{pmatrix}\begin{pmatrix}~V_1^*\\~V_2^*\end{pmatrix}~~{\rm and}~
 M_r=U_1\Sigma_1 V_1^*
\end{eqnarray*} 
  be SVDs
 of the matrices $M$ and its rank-$r$ truncation $M_r$, 
 respectively. [$\Sigma_2=O$ and $XY=M$ if 
 $\rank(M)=r$. The  $r$ columns of $V_1$ are the $r$ top right singular vectors of $M$.]
 Write 
   \begin{equation}\label{eqc12}
C_1=V^*_1H,~  
 C_2=V^*_2H.
 \end{equation}
Assume that $||H||_2 \le 1$ and $\rank(C_1) = r$. Then
\begin{equation}\label{eqerrnrm} 
 |||M-XY|||^2\le
 |||\Sigma_2|||^2+|||\Sigma_2C_2C_1^+|||^2.
 \end{equation}
 \end{theorem}

 
\begin{corollary}\label{copert0}
Under the assumptions of Thm. \ref{thpert1} and for
$\tau_{r+1}(M)$ of Thm.  \ref{thtrnc} it holds that
\begin{equation}\label{eqmmr1}
|||M-XY|||/\tau_{r+1}(M) \le (1+|||C_1^+|||^2)^{1/2} ~{\rm for}~C_1=V_1^*H.
\end{equation}

\end{corollary}

\begin{proof}
The corollary follows from (\ref{eqerrnrm}) because
$$|||\Sigma_2|||=\tau_{r+1}(M),~ |||C_2|||\le 1,~
{\rm and}~|||\Sigma_2C_2C_1^+|||\le |||\Sigma_2|||~~|||C_2|||~~|||C_1^+|||.$$  
\end{proof}

(\ref{eqmmr1}) implies
that the output {\em LRA} is optimal under both spectral and Frobenius matrix norms up to a factor of $(1+|||C_1^+|||^2)^{1/2}$.


    
\subsection{Impact of pre-multiplication on the errors of  {\em LRA} }\label{simppre}

 
 The following  theorem shows that 
the overall  error bounds of Alg. \ref{alg01}  
are dominated by the product of the norm $|||W|||$ and the
error norm bound  of Alg. \ref{alg1}.

\begin{theorem}\label{thst3}  {\rm See \cite{TYUC17}.} 
Let Alg. \ref{alg01} output  a matrix $XY$ for 
$Y=(FX)^+FM$ and let $m\ge k\ge l=\rank(X)$. Then
\begin{equation}\label{equvmpre}
 M-XY=W(M-XX^+M)~{\rm for}~W=I_m-X(FX)^+F,
\end{equation}
\begin{equation}\label{eqf} 
 |||M-XY|||\le |||W|||~~|||M-XX^+M|||,~~|||W|||\le  
 |||I_m|||+|||X|||~~|||F|||~~|||(FX)^+|||.
\end{equation} 
 \end{theorem}
 \begin{proof} 
 Recall that $Y=(FX)^+FM$ and notice that 
 $(FX)^+FX=I_l$
 if $k\ge l=\rank(FX)$. 
 Therefore,
 $Y=X^+M+(FX)^+F(M-XX^+M)$. Consequently, (\ref{equvmpre}) and  (\ref{eqf}) hold.
\end{proof} 
 
\begin{remark}\label{reprm}
  Deduce that
$|||W|||\le
|||I_m|||+|||F|||~|||F^+|||~|||X|||~|||X^+|||$ by combining bound (\ref{eqf})
with Lemma \ref{lehg}.
Hence $$|||W|||\le |||I_m|||+1$$ 
 if  both matrices $F$  and $X$ have orthonormal columns. The sketch matrix  $F$ is our choice, and for $X$ we ensure column  orthogonality   by properly choosing the matrix $T$ in Alg. \ref{alg01} 
 (see Remark \ref{re0}).
\end{remark}



 
%

\section{Output error norm bounds for dual  sketching algorithms}\label{serrranin}
 
 
Given the matrices 
$MHT^{-1}$ and $S^{-1}FM$,  
Alg. \ref{alg01} uses  
$O(kln)$ flops and hence
runs superfast for $kl\ll m$.
If also $l^2\ll m$ and $k^2\ll n$, then
for proper   Ultrasparse  matrices $F$ and $H$  we can  compute the matrices  $MHT^{-1}$ and $S^{-1}FM$ superfast  as well
and hence can perform entire Algorithm  \ref{alg01} superfast, involving much less than $mn$ entries of $M$ and other scalars,  already where $k\ll m$ and $l\ll n$.  Although such an {\em LRA}  fails for a worst case input,  it succeeds
whp under our dual model  in the case of any fixed, possibly sparse, well-conditioned  multipliers $F$ and $H$ of full rank
and a random input matrix $M$ that admits {\em LRA}. In view of the previous subsection we only need to prove this for Alg. \ref{alg1}.
  

\subsection{Auxiliary results}


\begin{lemma}\label{lemma:least_sing_bound}
Suppose that $G$ and $H$ are $r\times n$ and $n\times l$ matrices, respectively,  $r < l < n$, $GH$ has full rank $r$, and  $Q$ is an $r\times n$ matrix with orthonormal rows such that $Q$ and $G$ have  the same row space. Then 
\begin{eqnarray*}
\sigma_r(QH) \ge \frac{\sigma_r(GH)}{\sigma_1(G)}.
\end{eqnarray*}
\end{lemma}

\begin{proof}
Without loss of generality, suppose that $R\in \mathbb  R^{r\times r}$ and $G = RQ$. Then
\begin{eqnarray*}
\sigma_r(QH) = \sigma_r(R^{-1}GH) \ge \sigma_r(GH)\cdot \sigma_r(R^{-1}).
\end{eqnarray*}
Hence
 $\sigma_r(R^{-1}) = 1/\sigma_1(G)$
 because the matrices $R$ and $G$ share their singular values.
\end{proof}
\begin{lemma}\label{lemma:gauss_least_bound} 
Suppose that $G$ is an $r\times n$  Gaussian  
matrix, $H$ is an  $n\times l$ matrix with orthonormal columns,
 $n > l > r$, $l \ge 4$, 
 $Q$ is a matrix with orthonormal rows, and  $Q$ and $G$ share their row space. 
Fix two positive parameters $t_1$ and $t_2<1$. 
Then  
\begin{eqnarray*}\label{eqpr}
\sigma_r(QH) \ge  t_2 \cdot \frac{\sqrt{l}- \sqrt{r/l}+ \sqrt{1/l}}{t_1 + \sqrt{r} + \sqrt{n} }\cdot \frac{1}{e} \
\end{eqnarray*}
with a probability no less than $1 - \exp (-t_1^2/2) - (t_2)^{l-r}$.
\end{lemma}

\begin{proof}
The matrix  $GH$ has  distribution of an $r\times l$ Gaussian  
 matrix by virtue of 
 Lemma \ref{lepr3},
 and hence we can  assume that it has full rank (see Remark
\ref{refllrnk}).
Now recall Thm. \ref{thsignorm} and claim (ii) of
Thm. \ref{thsiguna} and obtain    
\begin{eqnarray*}
\prob \{ \sigma_1(G) > t_1 + \sqrt{r} + \sqrt{n}  \} < \exp(-t_1^2/2)~{\rm for}~t_1 \ge 0~{\rm and}
\end{eqnarray*}
\begin{eqnarray*}
\prob \Big\{ \sigma_r(GH) \le t_2\cdot\frac{l - r + 1}{e\sqrt{l}}  \Big\} \le (t_2)^{l - r}~{\rm for}~t_2 < 1~{\rm and}~l \ge 4. 
\end{eqnarray*}
Combine the latter  two inequalities, 
the union  
bound,
and Lemma \ref{lemma:least_sing_bound} and obtain  Lemma \ref{lemma:gauss_least_bound}.                                                                                                                                                                                                                                                                                                                                                                                                                                                                                                                                                                                                                                                                                                                                                                                                                                                                                                                                                                                                                                                                                              
\end{proof}

Lemma \ref{lemma:gauss_least_bound} implies that   $\sigma_r(QH)$ has at least order of $\sqrt{l/n}$ whp. 

\begin{corollary}\label{coro:least_sing_bound}
For $n, l, r, G, Q$, and  $H$ of Lemma \ref{lemma:gauss_least_bound},  let 
$n >  36r$  and  $ l > 22(r - 1)$. Then 
\begin{eqnarray*}
\prob \big\{  \sigma_r(QH) \ge \frac{1}{4} \sqrt{l/n}  \big\} > 1 - \exp\Big(-\frac{n}{72}\Big) - \exp\Big(-\frac{l-r}{20}\Big).
\end{eqnarray*} 
\end{corollary}
\begin{proof}  
Write $t_1: = \frac{1}{3}\sqrt{n} - \sqrt{r}$ and $t_2: = \frac{1}{3}~\frac{le}{l-r+1}$, recall that 
$n >  36r$  and  $ l > 22(r - 1)$, 
and then readily verify that $t_1 > \frac{\sqrt{n}}{6}$ and $\exp(0.05) > 0.95 > t_2 > 0$. Finally apply Lemma \ref{lemma:gauss_least_bound} under these bounds on $t_1$ and $t_2$.
\end{proof}

\begin{remark}\label{remark:least_sing_bound} 
We   
can extend
 our lower bounds of  Lemma \ref{lemma:gauss_least_bound} and Cor. \ref{coro:least_sing_bound} on 
$\sigma_r(QH)$   to the case 
of any matrix $H$  of full rank $l$ if we decrease these bounds
by a factor of $\sigma_l(H)$. 
\end{remark}


  
\subsection{Output errors of Alg. \ref{alg1}
for a matrix with random singular space}\label{sranssp}
 

\begin{assumption}\label{assump:mat_rand_space}
Let $r < n \le m$
(we can readily extend our study to the case where  $m < n$). Fix two constant matrices 
$$
\Sigma_r = 
\diag(\sigma_j)_{j=1}^r~{\rm and }~
\Sigma_{\perp} = 
\begin{bmatrix}
\diag(\sigma_j)_{j=r+1}^n\\
O_{m-n,n-r}
\end{bmatrix}\in \mathbb R^{(m-r)\times (n-r)}
$$
such that $\sigma_1\ge \sigma_2\ge \dots \ge \sigma_n \ge 0$, and $\sigma_r > 0$.
Let $G$ be an $r\times n$ Gaussian 
matrix 
and 
let $Q\in\mathbb R^{n\times r}$ and $Q_{\perp}\in\mathbb R^{(n-r)\times r}$ be  two matrices whose column sets 
make up orthonormal bases of the row space of $G$ and its orthogonal complement, respectively.
Furthermore, let 
\begin{equation}\label{eqassmp}
M = U\cdot\begin{bmatrix}\Sigma_r & 0\\ 0 & \Sigma_\perp \end{bmatrix} \cdot\begin{bmatrix}Q^*\\Q_\perp^*\end{bmatrix}
\end{equation}
be SVD where the
right
 singular space of $M$ is  random
 and $U$ is a  matrix
 with orthonormal columns.
\end{assumption}

\begin{remark}
 The matrix $G$ does not  uniquely define  
the matrices $Q$, $Q_{\perp}$, and $U$
under Assumption \ref{assump:mat_rand_space}  
and hence does not  uniquely define the matrix $M$, but this is immaterial
for our analysis. 
\end{remark}

\begin{theorem}\label{eqtherr} {\rm [Errors of Alg. \ref{alg1} for an input with a random singular space.]}
Suppose that $G$ is an $r\times n$ Gaussian  
matrix,
 $H\in\mathbb R^{n\times l}$ is a constant matrix, $n > 36r$, $l>22(r-1)$,   $r \le l < \min (m, n)$,
and Alg. \ref{alg1} applied  to the matrix $M$ of (\ref{eqassmp}) outputs  two  matrices $X$ and $Y$.

(i) If   $H$ has orthonormal columns,
then 
\begin{eqnarray*}
|||M - XY |||/\tilde \sigma_{r+1}(M) \le \sqrt{1 + 16n/l}
\end{eqnarray*}
with a probability no less than
$1 - \exp(-\frac{n}{72}) - \exp(-\frac{l-r}{20})$.

(ii) If   $H$ has full rank $l$, then
\begin{eqnarray*}
|||M - XY |||/\tilde \sigma_{r+1}(M) \le \sqrt{1 + 16\kappa^2(H)n/l}
\end{eqnarray*}
 with a probability no less than
$1 - \exp(-\frac{n}{72}) - \exp(-\frac{l-r}{20})$.
%
%

\end{theorem}

\begin{proof} We can assume that  the matrices $G$ and $GH$ have full rank $r$ (see Thm. \ref{thrnd} and Remark \ref{refllrnk}).
Consider SVD 
\begin{eqnarray*}
M = U\cdot\begin{bmatrix}\Sigma_r & 0\\ 0 & \Sigma_\perp \end{bmatrix} \cdot\begin{bmatrix}Q^*\\Q_\perp^*\end{bmatrix},
\end{eqnarray*}
  write $C_1 := Q^*H$, apply  Cor. \ref{copert0}, and deduce that
$$
|||M - XY |||/\tilde\sigma_{r+1}(M) \le \sqrt{1 + (||C_1^+||_2)^2}.
$$
(i) Recall from Cor. \ref{coro:least_sing_bound} that
\begin{eqnarray}
\prob \big\{  \sigma_r(C_1) = \sigma_r(Q^*H) \ge \frac{1}{4} \sqrt{l/n}  \big\} > 1 - \exp\Big(-\frac{n}{72}\Big) - \exp\Big(-\frac{l-r}{20}\Big).\label{ineq:rand_sub}
\end{eqnarray} 

\noindent (ii) Let $H = U_H\Sigma_HV_H^*$ be a compact SVD. Then
$\sigma_r(C_1)\ge \sigma_r(Q^*U_H)\sigma_r(H)$. 

Similarly to (\ref{ineq:rand_sub}) obtain that
\begin{align*}
\prob \big\{  \sigma_r(C_1) \ge \frac{1}{4} \sqrt{l/n} \cdot\sigma_r(H)  \big\} &\ge \prob \big\{  \sigma_r(Q^*U_H) \ge \frac{1}{4} \sqrt{l/n}\big\}\\
&> 1 - \exp\Big(-\frac{n}{72}\Big) - \exp\Big(-\frac{l-r}{20}\Big).
\end{align*}
Combine Thm.  \ref{thpert1}, equation $||Q_\perp||_2 = 1$, and the bound $\sigma_r(C_1) \ge \frac{1}{4} \sqrt{l/n} \cdot\sigma_r(H)$ and obtain
\begin{eqnarray*}
|||M - XY |||/\tilde \sigma_{r+1}(M) \le \sqrt{1 + ||Q_\perp^*HC_1^+ ||_2^2} \le \sqrt{1 + 16\kappa^2(H)n/l}.
\end{eqnarray*}
\end{proof}

Bound the output errors  of Algs. \ref{alg01} and 3.4
 by combining the estimates of this section and Sec. \ref{simppre} and by transposing the input matrix $M$.
 
  
 \subsection{Output errors of Range Finder for a perturbed factor-Gaussian input}\label{serrrang}
 

\begin{assumption}\label{assmp2} 
For  an $r\times n$ Gaussian  
matrix $G$ and a constant matrix
$A\in \mathbb R^{m\times r}$   of full rank $r< \min (m, n)$ define the matrices
 $B: = \frac{1}{\sqrt{n}}\cdot G$ 
 and  $\Tilde{M}: = AB$ and 
  call $M: = \Tilde{M} + E$ a perturbed  
   right factor-Gaussian 
   matrix if the
Frobenius norm of a perturbation matrix $E$ is   sufficiently small in comparison to $\sigma_r(A)$, as we specify in the next theorem.   
\end{assumption}


\begin{theorem}\label{therrfctr} {\rm [Errors of Range Finder for a  perturbed right factor-Gaussian matrix.]} 
Suppose that we are given an $r\times n$  Gaussian   matrix $G$ and constant matrices
$H\in\mathbb R^{n\times l}$, $A\in\mathbb R^{m\times r}$, and $E \in\mathbb R^{m\times n}$ for $r\le l < \min (m ,n)$.
Let $\tilde M$ be a right factor-Gaussian  matrix
 of Assumption \ref{assmp2}
and 
let $M = \Tilde{M} + E$  for a perturbation matrix $E$.
Apply Alg. \ref{alg1} to the  matrix $M$ with a test matrix $H$, 
and let $X$ and $Y$ denote the output matrices. Assume that $n > 36r$ and $l>22(r-1)$.

(i) Suppose that the matrix $H$ has orthonormal columns 
and  that $||E||_F \le \frac{\sigma_r(A)}{48\sqrt{n/l} + 6}$. Then
$$
|||M-XY|||/\tilde \sigma_{r+1}(M) \le \sqrt{1+ 100n/l}
$$
with a probability no less than $1 - \exp(-\frac{l-r}{20}) - \exp(-\frac{n-r}{20}) - \exp(-\frac{n}{72})$.
  
(ii) Assume  that $H$ has full rank and that  $||E||_F \le \frac{\sigma_r(A)}{12}\min(1,  \frac{1}{4\sqrt{n/l}\cdot\sigma_l(H) + 0.5})$
and let $\kappa(H) = ||H||_2||H^+||_2$ denote the spectral condition number of $H$. Then
$$
|||M-XY|||/\tilde \sigma_{r+1}(M)\le \sqrt{1 + 100\kappa^2(H)n/l } 
$$ 
with a probability no less than $1 - \exp(-\frac{l-r}{20}) - \exp(-\frac{n-r}{20}) - \exp(-\frac{n}{72}).$

\end{theorem} 


%



\begin{proof}

Recall that  the matrices $B$, $AB$, and $BH$  have full rank  with probability 1  
(see Thm. \ref{thrnd}) and assume that they do have full rank (see Remark \ref{refllrnk}). 
Let
\begin{eqnarray*}
M = \begin{bmatrix}U_r & U_{\perp} \end{bmatrix} \begin{bmatrix}\Sigma_r & 0 \\ 0 & \Sigma_\perp \end{bmatrix} \begin{bmatrix}V_r^T \\ V_\perp^T\end{bmatrix} 
\textrm{ and }
\tilde M = \begin{bmatrix}\tilde U_r & \tilde U_{\perp} \end{bmatrix} \begin{bmatrix}\tilde \Sigma_r & 0 \\ 0 & 0 \end{bmatrix} \begin{bmatrix}\tilde V_r^T \\ \tilde V_\perp^T\end{bmatrix} 
\end{eqnarray*}
be SVDs, where
$V_r$ and $\tilde V_r$ are the matrices of the $r$ top right singular vectors of $M$ and $\tilde M$, respectively.
Define $C_1 = V_r^TH$ and $\tilde C_1 := \tilde V_r^TH$ as in (\ref{eqc12}). 
Now Thm. \ref{thpert1}, that is, \cite[Thm. 9.1]{HMT11},
implies that
\begin{eqnarray*}
||| M - XY ||| /\tilde \sigma_{r+1}(M) \le \sqrt{1 + ||V_\perp^THC_1^+ ||_2^2}\label{ineq:fact_gauss_pert}.
\end{eqnarray*}

 Next we prove {\bf claim (i)}. Since $||V_\perp||_2 = ||H||_2 = 1$,  we only need to  obtain $\sigma_r(C_1) > \frac{1}{10}\sqrt{l/n}$  whp.
  
Recall that the columns of $V_r$ span the row space of an $r\times n$ Gaussian  
matrix $G$ and deduce from 
 Corollary \ref{coro:least_sing_bound}  that $\sigma_r(\tilde C_1) > \frac{1}{4}\sqrt{l/n}$ with a probability no less than $1 - \exp(-\frac{l-r}{20}) - \exp(-n/72)$.
It remains to estimate how   the perturbation $E$  alters the leading right  singular space of $\tilde M$  and how close is $\sigma_r(C_1)$  to $\sigma_r(\tilde C_1)$.
 
  If the norm of the perturbation matrix $||E||_F$ is sufficiently small, then
  by virtue of Lemma \ref{lemma:pert_sing_space} there exists a matrix $P$ such that  $\tilde V_r + \tilde V_\perp P$ and $V_r$ have the same column space and that, furthermore, $||P||_F \le \frac{2||E||_F}{\sigma_r(\tilde M) - 2||E||_F }$. This implies a desired bound on the differences of the smallest positive singular values of $C_1$ and $\tilde C_1$; next we  supply the details.
 
Claim (ii) of Thm. \ref{thsiguna} 
implies that  {\em whp} the $r$-th singular value of $\tilde M = AB$ is not much less than the $r$-th singular value of $A$.
Readily  deduce
from Cor. \ref{coro:least_sing_bound} that
\begin{eqnarray*} 
\prob \{ \sigma_r(\tilde M) < \sigma_r(A)/3 \} \le \prob \{\sigma_r(B) = \frac{1}{\sqrt{n}} \sigma_r(G) < 1/3\} \le e^{-(n-r)/20}.
\end{eqnarray*}
Hence $||E||_F \le \frac{\sigma_r(A)}{48\sqrt{n/l} + 6}\le \frac{\sigma_r (\tilde M) }{16\sqrt{n/l} + 2}$ with a probability no less than $1 - e^{-(n-r)/20}$, and so  $||P||_2 \le \frac{1}{8}\sqrt{l/n}$ for some  matrix $P$ of Lemma \ref{lemma:pert_sing_space} such that  $||P||_2 \le \frac{1}{8}\sqrt{l/n}$.  

Now let this holds, let  $\sigma_r(\tilde C_1) > \frac{1}{4}\sqrt{l/n}$, and deduce that
\begin{align}
\sigma_r(V_r^TH) &= \sigma_r\big((I_r + P^TP)^{-1/2}(\tilde V_r^T + P^T\tilde V_\perp^T)H\big) \label{eq:pf_of_thm} \\ 
&\ge \sigma_r\big((I_r + P^TP)^{-1/2}\big) \sigma_r(\tilde V_r^TH + P^T\tilde V_\perp^TH)\\
&\ge \frac{1}{\sqrt{1 + (\sigma_1(P))^2}} \big(  \sigma_r(\tilde C_1) - \sigma_1(P)\big) > \frac{1}{10} \sqrt{l/n}. \label{ineq:pf_of_thm2}
\end{align}
\noindent Equality (\ref{eq:pf_of_thm}) holds because the matrix $\tilde V_r + \tilde V_\perp P$ is normalized by $(I_r + P^TP)^{-1/2}$ (see Remark \ref{remark:pert_sing_space}) and has the same column span as $V_r$.
By applying the union bound deduce that inequality (\ref{ineq:pf_of_thm2}) holds with a probability no less than $1 - \exp(-\frac{l-r}{20}) - \exp(-\frac{n-r}{20}) - \exp(-\frac{n}{72})$.

To prove {\bf claim (ii)}, we essentially need to show that $\sigma_r(C_1) = \sigma_r(V_r^TH) \ge \frac{1}{10}\sqrt{l/n}\cdot \sigma_l(H)$, and then the claim will follow readily from inequality (\ref{ineq:fact_gauss_pert}).  
Let $H = U_H\Sigma_HV_H^T$ be compact SVD,
such that $U_H\in \mathbb R^{n\times l}$ , $\Sigma_H\in \mathbb R^{l\times l}$, and $V_H\in \mathbb R^{l\times l}$, and obtain that $\sigma_r(\tilde C_1) \ge \sigma_r(\tilde V_r^TU_H)\sigma_l(\Sigma_H)$ and 
\begin{eqnarray*}
\prob \{ \sigma_r(\tilde C_1) < \frac{1}{4}\sqrt{l/n}\cdot\sigma_l(H)\} < \exp(-\frac{l-r}{20}) + \exp(-\frac{n}{72}).
\end{eqnarray*}
Next  bound $\sigma_r(C_1)$ by showing that  the column spaces of $V_r$ and $\tilde V_r$ are sufficiently close 
to one another
if the perturbation $V_r-\tilde V_r$ is sufficiently small. 
Namely, assume that $||E||_F \le \frac{\sigma_r(A)}{12}$, and  then the assumptions of Lemma \ref{lemma:pert_sing_space} holds {\em whp}.
By applying the same argument as in the proof of claim (i), deduce that 
\begin{eqnarray*}
||E||_F \le \min\Big( \frac{\sigma_r (\tilde M)}{4}, \frac{\sigma_r (\tilde M) }{16\sqrt{n/l}\cdot\sigma_l(H) + 2}\Big) 
\end{eqnarray*}
with a probability no less than $1 - e^{-(n-r)/20}$.
It follows that  $||P||_2 \le \frac{1}{8}\sqrt{l/n}\cdot\sigma_l(H)$ for some matrix $P$ of Lemma \ref{lemma:pert_sing_space}. Hence $\sigma_r(C_1)\ge \frac{1}{10}\sqrt{l/n}\cdot \sigma_l(H)$ {\em whp}.
\end{proof}


 
\section{Numerical tests}\label{srndsmpl}

 
 Next we cover our tests  
 of dual superfast 
  variants of Algs. \ref{alg1}
  and \ref{alg01}
  specified in
 Remark \ref{re0}.
  As in Sec.  \ref{sexpr} we applied the  standard normal distribution function randn of MATLAB
  to generate Gaussian matrices. 
 We apply MATLAB function "svd()"    to calculate
 the $\epsilon$-rank  
for  $\epsilon=10^{-6}$ and perform the tests for Tables \ref{SuperfastTable}--\ref{tab8.10}  
on a 64-bit Windows machine with an Intel i5 dual-core 1.70 GHz processor by using custom programmed software in $C^{++}$ and compiled with LAPACK version 3.6.0 libraries, as in Sec. \ref{sexpr}.

  
\subsection{Input matrices}\label{ststmtrcs}

   
 We generated the following classes of 
 input matrices $M$ for testing {\em LRA} algorithms.
\medskip

{\bf Class I:} $M=U_M\Sigma_M V_M^*$,
where $U_M$ and $V_M$ are the Q factors
of the thin QR orthogonalization of  
$n\times n$ Gaussian matrices, 
$\Sigma_M=\diag(\sigma_j)_{j=1}^n$; 
 $\sigma_j=1/j,~j=1,\dots,r$,
$\sigma_j=10^{-10},~j=r+1,\dots,n$
(cf. [H02, Sec. 28.3]), 
 and  $n=256,
512,
1024$.
(Hence $||M||_2=1$ and 
$||M^+||_2=10^{10}$.)    

\medskip

{\bf Class II:}   
(i) The matrices $M$ of   the discretized single-layer Laplacian operator of  \cite[Sec. 7.1]{HMT11}:
$[S\sigma](x) = c\int_{\Gamma_1}\log{|x-y|}\sigma(y)dy,x\in\Gamma_2$,
for two circles $\Gamma_1 = C(0,1)$ and $\Gamma_2 = C(0,2)$  on the complex plane.
We arrived at the matrices      $M=(m_{ij})_{i,j=1}^n$,  
$m_{i,j} = c\int_{\Gamma_{1,j}}\log|2\omega^i-y|dy$ 
for a constant $c$,  $||M||=1$ and
 the arc $\Gamma_{1,j}$  of  $\Gamma_1$ defined by
the angles in the range $[\frac{2j\pi}{n},\frac{2(j+1)\pi}{n}]$.

(ii) The matrices  that approximate the inverse of a large sparse
matrix obtained from a finite-difference operator
of  \cite[Sec. 7.2]{HMT11}.

\medskip

{\bf Class III:} The dense  matrices of five families 
 from Class II of Sec. \ref{ststmtrcs0}. 

\medskip

We used $1024\times 1024$ SVD-generated input matrices of class I having numerical rank $r = 32$,  $400 \times 400$ Laplacian input matrices of class II(i)
having numerical rank $r = 36$,
 $408 \times 800$ matrices having numerical rank $r = 145$  and
representing finite-difference inputs  of class II(ii),
and $1000 \times 1000$ matrices of class III, having numerical rank 4, 6, 10, 12, and 25
for the matrices of the five  families {\em wing, baart, foxgood, shaw}, and {\em gravity}, respectively.




\subsection{Five families of Ultrasparse  multipliers $H$}\label{s17m} 

We generated our  $n\times (r+p)$  sketch matrices $H$ for random $p=1,2, \dots, 21$ by using
{\em  3-ASPH,
  3-APH (see the end of  Appendix \ref{ssrht}), and
 Random permutation matrices.}
  When  the overestimation parameter $p$ was considerable,
we actually computed {\em LRA} of numerical  rank larger than $r$, and so {\em LRA} was frequently  
closer to an input matrix than 
the optimal 
rank-$r$ approximation. Accordingly,
 the output error norms in our tests ranged from about $10^{-4}$ to $10^{4}$ {\em  relative to the optimal errors}.

We  obtained every 3-APH and every 3-ASPH
matrix by applying three Hadamard's recursive steps
(\ref{eqrfd}) followed by random column permutation defined by random permutation of the integers from 1 to $n$  inclusive. While generating a 3-ASPH matrix we also applied random scaling  with a diagonal matrix $D=\diag(d_i)_{i=1}^n$ where we have 
chosen the values of
independent identically distributed  {\em (iid)}  random variables $d_i$ sampled under the uniform
probability distribution  from the set
$\{-4, -3, -2, -1, 0, 
1 ,2, 3, 4\}$.
 
We used the following families of sketch matrices $H$:
(0)	Gaussian (for control),
(1)	sum of a 3-ASPH  and a permutation matrix,
(2) sum of a 3-ASPH  and two permutation matrices,
(3)	sum of a 3-ASPH  and three permutation matrices,
(4)	sum of a 3-APH  and three  permutation matrices, and
(5)	sum of a 3-APH  and two permutation matrices.

  
\subsection{Test results}\label{ststrslts}

 Tables
 \ref{SuperfastTable}--\ref{tab8.13} 
 display the average relative error norm  $ \frac{\|M - \tilde M \|_2}{\|M - M_{nrank}\|_2}$ in  our tests repeated 100 times for each class of input matrices and
 each 
size of an input matrix and multiplier $H$ for Alg.  \ref{alg1} or    
 for each size of an input matrix and a pair of left-hand and right-hand multipliers $F$ and $H$  for Alg.  \ref{alg01}.
 

In all our tests we applied the multipliers of the six  families of the previous subsection. 

   
 Tables  \ref{SuperfastTable}--\ref{tab8.10} display the average relative error norm for the output of
 Alg. \ref{alg1}; in  our tests 
 it ranged from about $10^{-3}$ to $10^{1}$.
 The numbers in parentheses in the first line of Tables \ref{tab8.9}
 and  \ref{tab8.10} show the numerical rank of input matrices.

 Tables  \ref{tab8.11}--\ref{tab8.13} display the average relative error norm for the output of
 Alg. \ref{alg01} 
applied  to the same input matrices from classes I--III as in our experiments for Alg. \ref{alg1}.

In these tests we used 
   $n\times \ell$ and $\ell\times m$ multipliers for $\ell = r+p$ and $k=c\ell$ for $c=1,2,3$ and random $p=1,2, \dots, 21$.

\begin{table}[htb]
\begin{center} 
\begin{tabular}{|c|c|c|c|c|c|c|}
\hline
			& \multicolumn{2}{|c|}{SVD-generated Matrices} & \multicolumn{2}{|c|}{Laplacian Matrices} & \multicolumn{2}{|c|}{Finite Difference Matrices}\\\hline
 \text{Family No.} & \text{Mean} & \text{Std} & \text{Mean} & \text{Std} & \text{Mean} & \text{Std} \\\hline
Family 0 & 4.52e+01 & 5.94e+01 & 6.81e-01 & 1.23e+00 & 2.23e+00 & 2.87e+00\\\hline
Family 1 & 3.72e+01 & 4.59e+01 & 1.33e+00 & 2.04e+00 & 8.22e+00 & 1.10e+01\\\hline
Family 2 & 5.33e+01 & 6.83e+01 & 1.02e+00 & 2.02e+00 & 4.92e+00 & 4.76e+00\\\hline
Family 3 & 4.82e+01 & 4.36e+01 & 7.56e-01 & 1.47e+00 & 4.82e+00 & 5.73e+00\\\hline
Family 4 & 4.68e+01 & 6.65e+01 & 7.85e-01 & 1.31e+00 & 3.53e+00 & 3.68e+00\\\hline
Family 5 & 5.45e+01 & 6.23e+01 & 1.03e+00 & 1.78e+00 & 2.58e+00 & 3.73e+00\\\hline
\end{tabular}
\caption{Relative error norms in tests for matrices of classes I and II}
 \label{SuperfastTable} 
\end{center}
\end{table}



\begin{table}
\begin{center}
\begin{tabular}{|c|c|c|c|c|c|c|}
\hline
			& \multicolumn{2}{|c|}{wing (4)} & \multicolumn{2}{|c|}{baart (6)} 
\\\hline
  
 \text{Family No.} & \text{Mean} & \text{Std} & \text{Mean}  & \text{Std}   \\\hline
 Family 0 	& 1.07e-03 & 6.58e-03 & 2.17e-02 & 1.61e-01 
\\\hline
 Family 1 	& 3.54e-03 & 1.39e-02  & 1.37e-02 & 6.97e-02 
\\\hline
 Family 2 	& 4.74e-03 & 2.66e-02  & 1.99e-02 & 8.47e-02  
\\\hline
 Family 3 	& 1.07e-03 & 5.69e-03  & 1.85e-02 & 8.74e-02  
\\\hline
 Family 4 	& 4.29e-03 & 1.78e-02  & 8.58e-03 & 5.61e-02  
\\\hline
 Family 5 	& 1.71e-03 & 1.23e-02  & 3.66e-03 & 2.38e-02 
\\\hline
\end{tabular}
\caption{Relative error norms for  input matrices of class III (of San Jose University database)}\label{tab8.9}
\end{center}
\end{table}

\begin{table}[htb] 
\begin{center}
\begin{tabular}{|c|c|c|c|c|c|c|}
\hline
			& \multicolumn{2}{|c|}{foxgood (10)} & \multicolumn{2}{|c|}{shaw (12)} & \multicolumn{2}{|c|}{gravity (25)}\\\hline
 
 \text{Family No.} & \text{Mean} & \text{Std} & \text{Mean}  & \text{Std} & \text{Mean} 	& \text{Std} \\\hline
 Family 0 	& 1.78e-01 & 4.43e-01  & 4.07e-02 & 1.84e-01  & 5.26e-01 & 1.24e+00\\\hline
 Family 1 	& 1.63e+00 & 3.43e+00  & 8.68e-02 & 3.95e-01  & 3.00e-01 & 7.64e-01\\\hline
 Family 2 	& 1.97e+00 & 4.15e+00  & 7.91e-02 & 4.24e-01  & 1.90e-01 & 5.25e-01\\\hline
 Family 3 	& 1.10e+00 & 2.25e+00  & 4.50e-02 & 2.21e-01  & 3.63e-01 & 1.15e+00\\\hline
 Family 4 	& 1.23e+00 & 2.11e+00  & 1.21e-01 & 5.44e-01  & 2.36e-01 & 5.65e-01\\\hline
 Family 5 	& 1.08e+00 & 2.32e+00  & 1.31e-01 & 5.42e-01  & 2.66e-01 & 8.22e-01\\\hline
\end{tabular}
\caption{Relative error norms for input matrices of class III (of San Jose University database)} \label{tab8.10}
\end{center}
\end{table}


\begin{table}[htb]
\begin{center}
\begin{tabular}{|c|c|c|c|c|c|c|c|}
\hline
		&	& \multicolumn{2}{|c|}{SVD-generated Matrices} & \multicolumn{2}{|c|}{Laplacian Matrices} & \multicolumn{2}{|c|}{Finite Difference Matrices}\\\hline	
 
\text{$k$} & \text{Class No.} & \text{Mean} & \text{Std} & \text{Mean}  & \text{Std} & \text{Mean} 	& \text{Std} \\\hline

\multirow{6}*{$\ell$} 
& Family 0 & 2.43e+03 & 1.19e+04 & 1.28e+01 & 2.75e+01 & 9.67e+01 & 1.48e+02 \\ \cline{2-8}
& Family 1 & 1.45e+04 & 9.00e+04 & 8.52e+03 & 8.48e+04 & 7.26e+03 & 2.47e+04 \\ \cline{2-8}
& Family 2 & 4.66e+03 & 2.33e+04 & 3.08e+01 & 4.07e+01 & 3.80e+02 & 1.16e+03 \\ \cline{2-8}
& Family 3 & 2.82e+03 & 9.47e+03 & 2.42e+01 & 3.21e+01 & 1.90e+02 & 3.90e+02 \\ \cline{2-8}
& Family 4 & 3.15e+03 & 7.34e+03 & 2.71e+01 & 4.69e+01 & 1.83e+02 & 2.92e+02 \\ \cline{2-8}
& Family 5 & 2.40e+03 & 6.76e+03 & 2.01e+01 & 3.56e+01 & 2.31e+02 & 5.33e+02  \\\hline

\multirow{6}*{$2\ell$}
& Family 0 & 5.87e+01 & 5.59e+01 & 7.51e-01 & 1.33e+00 & 3.17e+00 & 3.89e+00 \\ \cline{2-8}
& Family 1 & 7.91e+01 & 9.86e+01 & 3.57e+00 & 7.07e+00 & 1.55e+01 & 2.39e+01 \\ \cline{2-8}
& Family 2 & 5.63e+01 & 3.93e+01 & 3.14e+00 & 4.50e+00 & 5.25e+00 & 5.93e+00 \\ \cline{2-8}
& Family 3 & 7.58e+01 & 8.58e+01 & 2.84e+00 & 3.95e+00 & 4.91e+00 & 6.03e+00 \\ \cline{2-8}
& Family 4 & 6.24e+01 & 4.54e+01 & 1.99e+00 & 2.93e+00 & 3.64e+00 & 4.49e+00 \\ \cline{2-8}
& Family 5 & 6.41e+01 & 6.12e+01 & 2.65e+00 & 3.13e+00 & 3.72e+00 & 4.54e+00  \\\hline

\multirow{6}*{$3\ell$}
& Family 0 & 9.29e+01 & 3.29e+02 & 8.33e-01 & 1.54e+00 & --- & --- \\ \cline{2-8}
& Family 1 & 5.58e+01 & 4.20e+01 & 3.09e+00 & 4.08e+00 & --- & --- \\ \cline{2-8}
& Family 2 & 5.11e+01 & 4.94e+01 & 1.70e+00 & 2.08e+00 & --- & --- \\ \cline{2-8}
& Family 3 & 6.70e+01 & 8.27e+01 & 2.35e+00 & 2.96e+00 & --- & --- \\ \cline{2-8}
& Family 4 & 5.36e+01 & 5.74e+01 & 2.14e+00 & 3.76e+00 & --- & --- \\ \cline{2-8}
& Family 5 & 4.79e+01 & 4.58e+01 & 1.81e+00 & 2.94e+00 & --- & --- \\\hline

\end{tabular}
\caption{Relative error norms in tests for matrices of classes I and II }
\label{tab8.11}
\end{center}
\end{table}

\begin{table}[htb] 
\begin{center}

\begin{tabular}{|c|c|c|c|c|c|}
\hline
		&	& \multicolumn{2}{|c|}{wing (4)} & \multicolumn{2}{|c|}{baart (6)} \\\hline
 
\text{$k$} & \text{Class No.} & \text{Mean} & \text{Std} & \text{Mean}  & \text{Std}  \\\hline
\multirow{6}*{$\ell$}
& Family 0 & 1.70e-03 & 9.77e-03 & 4.55e+00 & 4.47e+01  \\ \cline{2-6}
& Family 1 & 3.58e+02 & 3.58e+03 & 1.42e-01 & 9.20e-01  \\ \cline{2-6}
& Family 2 & 2.16e-01 & 2.10e+00 & 1.10e-02 & 6.03e-02  \\ \cline{2-6}
& Family 3 & 7.98e-04 & 7.22e-03 & 4.14e-03 & 3.41e-02  \\ \cline{2-6}
& Family 4 & 5.29e-03 & 3.57e-02 & 2.22e+01 & 2.21e+02  \\ \cline{2-6}
& Family 5 & 6.11e-02 & 5.65e-01 & 3.33e-02 & 1.30e-01  \\\hline

\multirow{6}*{$2\ell$}
& Family 0 & 7.49e-04 & 5.09e-03 & 5.34e-02 & 2.19e-01  \\ \cline{2-6}
& Family 1 & 4.74e-03 & 2.32e-02 & 2.14e-02 & 1.33e-01  \\ \cline{2-6}
& Family 2 & 3.01e-02 & 2.34e-01 & 1.26e-01 & 7.86e-01  \\ \cline{2-6}
& Family 3 & 2.25e-03 & 1.38e-02 & 5.91e-03 & 2.63e-02  \\ \cline{2-6}
& Family 4 & 3.94e-03 & 2.54e-02 & 1.57e-02 & 6.71e-02  \\ \cline{2-6}
& Family 5 & 2.95e-03 & 1.47e-02 & 1.58e-02 & 1.20e-01  \\\hline

\multirow{6}*{$3\ell$}
& Family 0 & 4.59e-03 & 2.35e-02 & 1.50e-02 & 7.09e-02  \\ \cline{2-6}
& Family 1 & 5.96e-03 & 2.82e-02 & 7.57e-03 & 4.84e-02  \\ \cline{2-6}
& Family 2 & 1.74e-02 & 1.06e-01 & 6.69e-03 & 2.97e-02  \\ \cline{2-6}
& Family 3 & 3.07e-03 & 3.07e-02 & 1.16e-02 & 5.16e-02  \\ \cline{2-6}
& Family 4 & 2.57e-03 & 1.47e-02 & 2.35e-02 & 9.70e-02  \\ \cline{2-6}
& Family 5 & 4.32e-03 & 2.70e-02 & 1.36e-02 & 5.73e-02  \\\hline


\end{tabular}
\caption{Relative error norms for  input matrices of class III (of San Jose University database) }
\label{tab8.12}
\end{center}

\end{table}

\begin{table}[htb]
\begin{center}
\begin{tabular}{|c|c|c|c|c|c|c|c|}
\hline
& & \multicolumn{2}{|c|}{foxgood (10)} & \multicolumn{2}{|c|}{shaw (12)} & \multicolumn{2}{|c|}{gravity (25)}\\\hline
 
\text{$k$} & \text{Class No.} & \text{Mean} & \text{Std} & \text{Mean}  & \text{Std} & \text{Mean} 	& \text{Std} \\\hline
\multirow{6}*{$\ell$}
& Family 0 & 5.46e+00 & 1.95e+01 & 8.20e-01 & 4.83e+00 & 8.56e+00 & 3.33e+01 \\ \cline{2-8}
& Family 1 & 8.51e+03 & 1.88e+04 & 1.12e+00 & 5.75e+00 & 1.97e+01 & 1.00e+02 \\ \cline{2-8}
& Family 2 & 5.35e+03 & 1.58e+04 & 1.93e-01 & 1.51e+00 & 8.79e+00 & 4.96e+01 \\ \cline{2-8}
& Family 3 & 6.14e+03 & 1.74e+04 & 4.00e-01 & 1.90e+00 & 7.07e+00 & 2.45e+01 \\ \cline{2-8}
& Family 4 & 1.15e+04 & 2.33e+04 & 2.95e-01 & 1.71e+00 & 4.31e+01 & 3.80e+02 \\ \cline{2-8}
& Family 5 & 7.11e+03 & 1.87e+04 & 2.18e-01 & 9.61e-01 & 6.34e+00 & 2.59e+01\\\hline

\multirow{6}*{$2\ell$}
& Family 0 & 2.70e-01 & 7.03e-01 & 5.54e-02 & 2.62e-01 & 5.34e-01 & 1.59e+00 \\ \cline{2-8}
& Family 1 & 5.24e+02 & 5.19e+03 & 4.67e-02 & 2.35e-01 & 1.38e+01 & 1.31e+02 \\ \cline{2-8}
& Family 2 & 2.45e+00 & 3.47e+00  & 8.31e-02 & 6.37e-01 & 5.47e-01 & 1.69e+00 \\ \cline{2-8}
& Family 3 & 2.43e+00 & 3.74e+00  & 1.24e-01 & 8.52e-01 & 5.10e-01 & 1.24e+00 \\ \cline{2-8}
& Family 4 & 2.17e+00 & 2.92e+00  & 1.76e-01 & 8.76e-01 & 2.60e-01 & 7.38e-01 \\ \cline{2-8}
& Family 5 & 2.10e+00 & 3.34e+00 & 1.26e-01 & 5.99e-01 & 5.68e-01 & 1.46e+00 \\\hline

\multirow{6}*{$3\ell$}
& Family 0 & 2.62e-01 & 8.16e-01 & 4.49e-02 & 1.64e-01 & 4.59e-01 & 1.38e+00 \\ \cline{2-8}
& Family 1 & 2.72e+00 & 4.60e+00 & 6.84e-02 & 3.43e-01 & 3.44e-01 & 8.60e-01 \\ \cline{2-8}
& Family 2 & 2.42e+00 & 3.92e+00 & 8.26e-02 & 5.38e-01 & 6.89e-01 & 2.15e+00 \\ \cline{2-8}
& Family 3 & 3.22e+02 & 3.20e+03 & 6.06e-02 & 2.95e-01 & 5.26e-01 & 1.17e+00 \\ \cline{2-8}
& Family 4 & 1.91e+00 & 3.36e+00 & 6.61e-02 & 3.36e-01 & 6.19e-01 & 1.54e+00 \\ \cline{2-8}
& Family 5 & 2.73e+00 & 6.90e+00 & 5.72e-02 & 2.39e-01 & 7.22e-01 & 1.59e+00 \\\hline

\end{tabular}
\caption{Relative error norms for  input matrices of class III (of San Jose University database)}
\label{tab8.13}
\end{center}

\end{table}

  
\section{Conclusions}\label{scncl} 

Next we list some promising directions
for further  study of superfast
LRA.

\cite{PLa}
 proposed and tested 
  iterative refinement of LRA, and \cite{GKLPPZa}
has greatly advanced that study incorporating
 the superfast  algorithms that we studied here.
 
\cite{GPa} proposed, analyzed, and tested  superfast algorithms for  matrix norm estimation.

  Least squares approximation of multivariate functions (LSAMF) was reduced in 
\cite{ALS24}
 to 
the solution of 
 Linear Least Squares Problem (LLSP) and further to LRA. Namely,
 \cite{ALS24}  reduced LSAMF to LLSP, which we state as follows: for given $M\in \mathbb C^{m\times n}$ and ${\bf b}\in \mathbb C^{m\times 1}$ where $m>n$, compute or approximate
  a vector
\begin{equation}\label{eqllsp}  
{\bf z}={\rm argmin}_{\bf y}||M{\bf y}-{\bf b}||.
\end{equation}
Triangular inequality implies that
for an LRA 
 $\tilde M$ of $M$ we have 
$$ ||M{\bf y}-{\bf b}||\le ||\tilde M{\bf y}-{\bf b}||+||(\tilde M-M){\bf y}||.$$
 
\cite{ALS24}  showed that $||(\tilde M-M){\bf y}||=0$
in the case of LSAMF problem,  and so we can
 solve it superfast if  LRA can
be computed superfast, e.g., with our algorithms. 

Indeed, suppose that our superfast LRA  algorithm has computed  $P$ and $S$ such that $\tilde M=PS\in \mathbb R^{m\times n}$, 
$P\in \mathbb R^{m\times r}$, 
$S\in \mathbb R^{r\times n}$ and $r\ll \min\{m,n\}$. Then  compute a solution
$${\bf x}={\rm argmin}_{\bf v}||P{\bf v}-{\bf b}||$$ of an LLSP for a matrix $P$ of  smaller size
$m\times r$ and  immediately recover the solution
$\tilde {\bf z}$ of the LLSP
  $ {\bf z}={\rm argmin}_{\tilde {\bf y}}||\tilde M\tilde {\bf y}-{\bf b}||$ for $\tilde M=PS\in \mathbb R^{m\times n}$, 
$P\in \mathbb R^{m\times r}$, 
$S\in \mathbb R^{r\times n}$,
and $r\ll \min\{m,n\}$ 
from the
underdetermined 
linear system
${\bf x}=S{\bf z}$.
Finally, extend  this solution of LLSP to superfast solution of  LSAMF by following \cite{ALS24}.


\medskip
\medskip
\medskip

{\bf \Large Appendix} 
\appendix


\section{Small families of hard inputs for 
superfast LRA}\label{shrdin}

  Any  algorithm
  not involving all entries of an input matrix 
  fails on the following input families.

\begin{example}\label{exdlt} 
 Let  $\Delta_{i,j}$ denote an $m\times n$ matrix
 of rank 1  filled with 0s except for its $(i,j)$th entry filled with 1. The $mn$ such matrices $\{\Delta_{i,j}\}_{i,j=1}^{m,n}$ form a family of  $\delta$-{\em matrices}.
We also include the $m\times n$ null matrix $O_{m,n}$
filled with 0s  into this family.
Now fix any superfast algorithm; it does not access the $(i,j)$th  
entry of its input matrices  for some pair of $i$ and $j$. Therefore, it outputs the same approximation 
of the matrices $\Delta_{i,j}$ and $O_{m,n}$,
with an undetected  error at least 1/2.
Arrive at the same conclusion by applying the same argument to the
set of $mn+1$ small-norm perturbations of 
the matrices of the above family and to the                                                                                                                                   
 $mn+1$ sums  
of the latter matrices with  any
   fixed $m\times n$ matrix of low rank.
Finally, the same argument shows that 
a posteriori estimation of the output errors of an LRA algorithm
applied to these input families
cannot run superfast. 
\end{example}

The example  covers randomized LRA algorithms as well. Indeed, let an LRA algorithm miss 
 an entry of an input matrix with a  probability $p>0$.
Apply this algorithm to two matrices
of low rank whose difference at  this entry is equal to a large constant $C$.   
Then
with a  
probability $p$ the algorithm  
has error at least $C/2$  at this entry for a least one of these two matrices. 


\section{LRA with dense and Abridged SRHT random  multipliers}
\label{sspfcrd}
                                                                                                                                                                                                                                                                                                                                                                                                                                                                                                                                                                                                                                                                                                                                                                                                                                                                                                                                                                                                                                                                                                                                                                                                                                                                                                                                                                                                                                                                                                                                                                                                                                                                                                                                                                                                                            
\subsection{The generalized Nystr{\"o}m  LRA algorithm (briefly)}
\label{sdrnd}


 
The  {\em 
generalized Nystr{\"o}m}  algorithms of \cite{TYUC17,N20} specify Alg. \ref{alg01} for random sketching LRA.\footnote{\cite{N20}  modifies the algorithms of \cite{TYUC17}  to
 improve numerical stability of rank-$r$ approximation of  
a matrix $M\in \mathbb C^{m\times n}$
 for 
$r\le  n\le m$, but in our tests we departed from \cite{TYUC17}
and accelerated it by using
Abridging SRHT multipliers $F$ and $H$.}
  In the best studied case  the multipliers $F$ and $H$ are  Gaussian random
  matrices, filled with independent
 standard Gaussian (normal) random variables; then for  
 $k=2r+1$ 
 and $\ell=2k$  the expected Frobenius error norm  of the error matrix $M-XYZ$ is within a factor of 2 from optimal (not exceeding the  optimal one for rank-$\frac{k}{2}$ approximation).                                                                                                                                                                                                                                                                                                                                                                                                                                                                                                                                                                                                                                                                                                                                                                                                                                                                                                                                                                                                                                                                   
The algorithm  involves $O((k+\ell)mn)$ 
 flops, $mn$ entries of $M$, and about 
 $km+n\ell+k\ell$ other scalars,   overwriting some entries of $M$. 
 
With SRHT or SRFT (rather than Gaussian) multipliers $F$ and $H$, one obtains  LRA (with a little higher error probability) by using
$O(mn\log(k\ell)))$ flops, $mn$ entries of $M$, and 
 $mk\log(k)+n\ell\log(\ell)$ other
 scalars, also  overwriting some  entries of $M$.
 
Like the other random sketching LRA algorithms,
this one is superfast except for the stage of the computation of the sketches $FM$  and  $MH$ and can be readily modified to
bound the relative error of the output LRA by any  $\epsilon>0$, with  estimated sketches' size increasing fast as  $\epsilon\mapsto 0$.

                                                                                                                                                                                                                                                                                                                                                                                                                                                                                                                                                                                                                                                                                                                                                                                                                                                                                                    
\subsection{Abridged SRHT  matrices}\label{ssrht}

 
With  {\em sparse subspace embedding} \cite{CDDRa,CFSa},
\cite[Sec. 3.3]{TYUC19},
\cite[Sec. 9]{MT20} one obtains significant acceleration but still does not yield superfast algorithms. 
According to \cite{L09},
 such acceleration tends to
 make the accuracy of the output LRAs somewhat less reliable, although \cite{CFSa} partly overcomes
 this  problem for incoherent matrices. One can  multiply a  matrix by Subsampled Randomized Hadamard or Fourier Transform ({\em SRHT or SRFT)} dense matrices
 towards  incoherence \cite{CFSa}, but this step is not superfast.
We make it superfast by applying {\em Abridged SRHT multipliers}, 
specified below. 
Thm.  \ref{eqtherr}   provides semi-empirical   support for the output accuracy
of the resulting LRA.   \cite{PLSZb} 
covers earlier
formal and empirical study of Abridged SRHT matrices and their
application to LRA.

We  proceed
  by  means of abridging 
    the classical 
recursive processes of the generation  of 
SRHT  matrices,
  obtained from   the  $n\times n$ dense
matrices $H_n$ of  Walsh-Hadamard transform for $n=2^t$
(cf. \cite[Sec. 3.1]{M11}).
Namely, we
write 
 

\begin{equation}\label{eqrfd}
H_{d,0}:=I_{n/2^d},~
H_{d,i+1}:=\begin{pmatrix}
H_{d,i} & H_{d,i} \\
H_{d,i} & -H_{d,i}
  \end{pmatrix}
  ~{\rm for}~i=0,1,\dots,d-1. 
\end{equation}
The $n\times n$ matrix $H_{d,d}$ is orthonormal 
up to scaling and
 has $2^d$ nonzero entries 
in every row and  column. 
Now define  $l\times n$ 
 matrix $R$ of uniform random  sampling of $l$ out of $n$ columns, the  $n\times n$ diagonal matrix 
  $D$  whose $n$ diagonal entries are independent random signs, i.e., random variables uniformly distributed on the pair
 $\{\pm 1\}$, and the $d$-{\em Abridged
 subsampled randomized Hadamard transform (SRHT)}
 matrix  $\sqrt{2^d/l~}RH_{d,d}D$, with
  $ 2^d$
nonzero entries
$\pm  \sqrt{2^d/l}$   per column. We can multiply $M\in \mathbb C^{m\times n}$  by that matrix by
involving $m2^d$ 
entries of $M$
and $(2^d-1)m$ flops.

 Likewise, by   sampling $k$ random rows of $H_{d,d}$ and multiplying 
 the resulting matrix by $\pm  \sqrt{2^d/k}D$  obtain 
$d$-Abridged $k\times m$ SRHT matrix having
$2^d$ nonzero entries $\pm  \sqrt{2^d/k}$ per row. We can multiply it by $M$ by involving  $n2^d$ 
entries of $M$
and $(2^d-1)n$ flops, where a flop amounts to addition, subtraction,
or multiplications by 
$\sqrt{2^d/l}$ or $\sqrt{2^d/k}$.
In particular 
this means  $8n$  (resp. $8m$) entries of $M$ and $7n$
(resp. $7m$) flops in our numerical tests, where we always fix $d=3$.

A $t$-Abridged SRHT matrix   turns into
an SRHT matrix.

Table \ref{tabalgb1} displays the computational cost of randomized rank-$r$ approximation by means of the algorithm of \cite{N20}, when it uses Gaussian, SRHT and $d$-Abridged   
SRHT multipliers $F$ and $H$.


\begin{table}[h!]
\centering
\begin{tabular}
{|l|c|c|c|}
\toprule
Multipliers $F$ and $H$& Input entries involved & Other entries involved &  Flops involved  \\
\midrule
Gaussian & $mn$& $km+n\ell+kl$ & $O((k+l)mn)$ \\

SRHT& $mn$ & $O(mk\log(k)+nl\log(l))$ & $O(mn\log(kl))$\\
$d$-Abridged SRHT &$2^d(m+n)$ & $2^d(k+l)$ & $(2^d-1)(m+n)$ \\
\bottomrule
\end{tabular}
\caption{Computational cost of the {\em 
generalized Nystr{\"o}m} algorithm  of \cite{N20} with 
Gaussian, SRHT, and $d$-Abridged SRHT multipliers $F$ and $H$.}\label{tabalgb1}
\end{table}



We also refer to $d$-abridged SRHT matrices as  
$d$-Abridged Scaled and Permuted 
 Hadamard $d$-{\em ASPH} matrices and call them $d$-{\em APH} matrices if random scaling is  omitted, that is, if the multiplier $D$, filled with $\pm 1$,
 is replaced with the identity matrix, and then randomization
 only comes from random permutation matrix. 


\section{The spectral norms of a Gaussian  matrix and its pseudo inverse}\label{snrmg}

Hereafter
$\Gamma(x)=
\int_0^{\infty}\exp(-t)t^{x-1}dt$
denotes the Gamma function; 
 $\nu_{p, q}$ and $\nu^+_{p, q}$ denote the random variables 
representing the spectral norms of a $p\times q$ Gaussian random matrix and its Moore-Penrose pseudo inverse, respectively.


\begin{theorem}\label{thsignorm} {\rm [Spectral norms of a Gaussian matrix. 
 See  \cite[Thm. II.7]{DS01}.]}

  {\rm Probability}$\{\nu_{m,n}>t+\sqrt m+\sqrt n\}\le
\exp(-t^2/2)$ for $t\ge 0$, 
 $\mathbb E(\nu_{m,n})\le \sqrt m+\sqrt n$.
\end{theorem}


\begin{theorem}\label{thsiguna} 
 {\rm [Spectral norms of the pseudo inverse of a Gaussian matrix.]} 

(i)  {\rm Probability} $\{\nu_{m,n}^+\ge m/x^2\}<\frac{x^{m-n+1}}{\Gamma(m-n+2)}$
for $m\ge n\ge 2$ and all positive $x$,

(ii)  {\rm Probability} $\{\nu_{m,n}^+\ge  	 
  t\frac{e\sqrt{m}}{m-n+1}\}\le t^{n-m}$
  for all $t\ge 1$ provided that $m\ge 4$, 
  
(iii) $\mathbb E(\nu^+_{m,n})\le \frac{e\sqrt{m}}{m-n}$ 
provided that $m\ge n+2\ge 4$,
 
\end{theorem}


\begin{proof}
 See \cite[Proof of Lemma 4.1]{CD05} for claim (i),   
\cite[Prop.
 10.4 and Eqns. (10.3) and (10.4)]{HMT11} for claims (ii) and (iii), 
and \cite[Thm. 3.3]{SST06} for claim (iv).
\end{proof}
 

Thm. \ref{thsiguna}
implies reasonable probabilistic upper bounds on the norm 
 $\nu_{m,n}^+$,
even where the integer $|m-n|$ is close to 0;  whp the upper bounds of Thm. \ref{thsiguna}
on the norm $\nu^+_{m,n}$ decrease very fast as the difference $|m-n|$ grows from 1.


\section{Randomized  pre-processing of LRA}\label{srndoprpr}


The following simple results 
from \cite[Section 8.2]{PLSZa} 
 show that pre-processing with 
Gaussian multipliers $X$ and $Y$ transforms any matrix that admits $LRA$
into a perturbation of a factor-Gaussian matrix.
We write $A  \preceq B$ (resp. $A \succeq B$) to show that $A$ is  statistically 
less (resp.   
 greater) than or equal to $B$.

\begin{theorem}\label{thquasi}  
Consider five integers $k$, $l$, $m$, $n$, and  $r$
satisfying the bounds
$r\le k\le m,~r\le l\le n$,
an $m\times n$ well-conditioned matrix $M$ 
of rank $r$, $k\times m$ and $n\times l$ Gaussian matrices $G$ and $H$, respectively,
and the norms $\nu_{p,q}$
and $\nu_{p,q}^+$ of Appendix \ref{snrmg}.
Then 

(i) $GM$ is a left factor-Gaussian matrix of expected rank $r$ such that
$$||GM||_2\preceq ||M||_2~\nu_{k,r}~{\rm  
and}~||(GM)^+||_2\preceq ||M^+||_2~\nu_{k,r}^+,$$

(ii) $MH$ is a right factor-Gaussian matrix of expected  rank $r$ such that
$$||MH||_2 \preceq ||M||_2~\nu_{r,l}~{\rm 
and}~ 
||(MH)^+||_2\preceq ||M^+||_2~\nu_{r,l}^+,$$ 
(iii) $GMH$ is a two-sided 
factor-Gaussian matrix of expected  rank 
$r$
such that $$||GMH||_2\preceq ||M||_2~\nu_{k,r}\nu_{r,l}~{\rm 
and}~||(GMH)^+||_2\preceq ||M^+||_2~\nu_{k,r}^+\nu_{r,l}^+.$$
\end{theorem}
\begin{remark}\label{reprprp}
According to this theorem we can readily extend our results on $LRA$ of perturbed 
factor-Gaussian matrices to all matrices that admit $LRA$ and are pre-processed with Gaussian multipliers.
We can perform such pre-processing
 superfast if we replace  
 Gaussian multipliers
 with 3-abridged SRHT
 orthogonal multipliers;
this replacement has little affected the output accuracy of LRA for a large class of  matrices $M$
  in the tests in  \cite{GKLPPZa}.
\end{remark}


\section{Proof of Thm. \ref{thncl}}\label{sncl}

\begin{proof}  $\sigma_j(G)\le \sigma_j(M)$ for all $j$
 because $G$ is a submatrix of $M$. Hence 
  $\epsilon$-rank$(G)\le \epsilon$-rank$(M)$ for all nonnegative $\epsilon$,
   and  in particular $\rank(G) \le \rank(M)$.                                                                  
 
 Now let $M=M'=CUR$. Then clearly $$\rank(M)\le \rank(U)=\rank(G_{r}^+)=\rank(G_{r})\le \rank(G),$$
 Hence  $$\rank(G) \ge \rank(M),~{\rm and~so}~ 
 \rank (G)=\rank (M)~{\rm if}~M'=M.$$
 
It remains to deduce that $$M=CG_{r}^+R~{\rm if}~\rank(G) =\rank(M):=r,$$
but in this case  $G_{r}=G$, and so 
$$\rank(CG_{r}^+R)=\rank(C)=\rank(R)=r.$$
Hence the rank-$r$ matrices $M$ and $CG_{r}^+R$ share  their rank-$r$ submatrices 
$C$
and $R$.
\end{proof}

 \begin{remark}\label{rencl}
 Can we extend the theorem by proving that  $M'\approx M$
 if and only if   $\epsilon$-rank$(G) = \epsilon$-rank$(M)$ for a small positive $\epsilon$? The ``only if'' claim  cannot be extended, e.g., for
 $$
 M = 
 \begin{pmatrix}
 1 & 0 & 0\\
 0 & \epsilon & 0\\
 0 & 1 & 0
 \end{pmatrix},$$
$\epsilon\approx 0$, and the $2\times 2$ leading submatrix  $G$ of $M$.
Indeed, $\rank(M)=
 \rank(G)=2$, and so Thm. \ref{thncl} implies
 that $M'=M$, while
 $\epsilon$-$\rank(M)=2>\epsilon$-$\rank(G)=1.$
 \end{remark}


 
\medskip


\noindent {\bf Acknowledgements:}
Our work has been supported by NSF Grants 
CCF--1116736,
 CCF--1563942 and CCF--1733834
and PSC CUNY Awards  66720-00 54 and 68452-00 56.
We are also grateful to E. E. Tyrtyshnikov for the challenge
of formally supporting empirical power of C--A iterations,
to N. L. Zamarashkin for his comments on \cite{PLSZ17},   and to
  S. A. Goreinov, 
 I. V. Oseledets, A. Osinsky,  E. E. Tyrtyshnikov, and N. L. Zamarashkin for  reprints and pointers to relevant bibliography.



\end{document}